# New Generating and Counting Functions of Prime Numbers applied to approximate Euler´s Product, Chebyschev 2nd class function, Riemann Zeta function, and applying the Least Action Principle to find non-trivial roots of Riemann´s Zeta

**Eduardo Stella, Celso L. Ladera and Guillermo Donoso**

**Departamento de Física, Universidad Simón Bolívar**

**Valle de Sartenejas, Baruta**

**Caracas 1089 Venezuela**

*Corresponding author: clladera@usb.ve*

*Prime numbers seem to be predestined to become the eventual paradigm of codification of the Universe*

**Abstract** We introduce novel prime numbers functions including an exact Generating Function and a Discriminating Function of Prime Numbers, neither based on prime number tables nor on algorithms. Instead, these functions are defined in terms of ordinary elementary functions, therefore having the advantage of being readily calculable. Also presented are four applications of our new Prime Numbers Generating Function, namely: (i) obtaining a new formula for counting Prime Numbers, (ii) obtaining an approximant to Euler´s product formula, (iii) obtaining an approximant to Riemann Zeta function $\zeta(\sigma,\tau)$ based on our prime numbers Discriminating Function, (iii) an accurate approximant to Chebyshev function of second class in terms of our Primes Numbers Generating function, and the application of this approximant in well-known sharp estimates related to the validity of the Riemann Hypothesis. We also apply the variational calculus of classical mechanics to obtain non-trivial roots of our good accuracy approximated Riemann zeta function, obtained in this work, in an original and novel approach. To that end a variational test function based on the modulus squared of the approximated Riemann function is defined, and then Hamilton Principle is applied to find non-trivial roots of Riemann's zeta with great accuracy along the critical line defined in the Riemann Hypothesis, in a completely original way and. We optimize our analytical variational procedure by using a second, more general test function, that depends explicitly upon the abscissa variable $\sigma$, and present a second procedure to obtain non-trivial roots with good accuracy. Our method even allows us to define a symmetrical function $F(\sigma)$ that behaves analogously to the well-known Functional Equation of the Riemann zeta function.

## 1 Introduction

Probably the best arguments for studying prime numbers were written by the mathematician A. E. Ingham in his book *The Distribution of Prime Numbers* [1]: "The prime numbers derive their peculiar importance from the 'fundamental theorem of arithmetic' that a composite number can be expressed in one and only one way as a product of prime factors. A problem which presents itself at the very threshold of mathematics is the question of the distribution of the primes among the integers. Although the series of prime numbers exhibits great irregularities of detail, their distribution among integers is found to possess certain features of regularity, which could be formulated in precise terms thus becoming the subject of mathematical investigation." Certainly those regularities of prime numbers have been the subject of rigorous investigation by mathematicians in the last three hundred years, among them L. Euler, K. F. Gauss, G. L. Dirichlet, B. Riemann, C-J. de la Vallée-Poussin, and more recently by G. H. Hardy, A. Selberg, P. Erdös, E. Bombieri, plus other influential mathematicians of the last century [2, 3]. For instance in 1837 Dirichlet considered prime integers *a*, *q* such that $2 \leq a \leq q$, and proved that any arithmetic progression distribution:

$$A(a,q) = \{a,\ q+a,\ 2q+a, 3q+a, 4q+a, \dots\}, \qquad (1)$$

contains infinitely many primes; an example being the primes in the distribution
$A(2,5) = \{\mathbf{2}, \mathbf{7}, 12, \mathbf{17}, 22, 27, 32, \mathbf{37}, 42, \mathbf{47}, \dots\}$.

A second example of attempting to find a prime numbers distribution is the *Bertrand Conjecture*, later a theorem proved by Chebyschev [4], that states: *There is at least one prime number between n and 2n for all $n \in \mathbb{N}$*. Important efforts to obtain approximate mathematical formulae for functions of prime numbers persist to present days [5-7], and one of the objectives of the present work is to obtain a new Prime Numbers Generating function that should be analytic, opening a door for new applications to the rest of the prime number based functions. Today some modern encryption systems –not a subject of this work – exploit that "fundamental theorem of arithmetic" mentioned by Ingham, namely: to codify information in terms of the prime factors of a large composite number that it is either hard to decompose, or very difficult to prove to be a prime, e.g. a well-known case was the factorization of the RSA 129-digit number proposed by M. Gardner in 1977 [8], which took several years to be solved. Prime numbers and their functions are said to belong to an advanced world of intellectual conceptions, certainly to Analytical Number Theory, Cryptography, and to be related to some other subfields of natural and formal sciences, where these numbers and their functions could play some role e.g. in Discrete Classical Mechanics, Quantum Mechanics, Statistical Physics, Chaos and Fractals, or in Nuclear Physics [9, 10].

In Section 2 of this work we introduce and exploit a new mathematical function, whose purpose is generating a prime numbers distribution. Defined in terms of elementary functions this Prime Numbers Generating function has the advantage of being defined in terms of ordinary functions, enabling us to apply both, algebraic procedures and tools of differential and integral calculus when using it, in an altogether original way; it of course enables us to operate and calculate within the entirely Discrete Classical Mechanics analytic formalism. Our Prime Numbers Generating function is not only validated here in finite intervals, but it is actually used in the construction of approximations to other emblematic functions of Number Theory, such as an exact Prime Numbers Counting function. In Section 3 we also present an approximant to Euler's product formula. It is an auto-consistent relation that, being based on our Prime Numbers functions, does not require the use of prime numbers tables or algorithms. In Section 4 we introduce our analytic approximant to Riemann zeta function of complex numbers $s=(\sigma,\tau)$, for the abscissa values $\sigma>0$, and particularly to approximate its important Analytic Continuation into that particular region [2, 4], and to obtain the non-trivial roots of that function along finite domains of the so-called critical line. Section 5 is devoted to apply Hamilton Principle of physics to obtain non-trivial roots of Riemann Zeta function by constructing a Lagrangian function based on a test function which involves the modulus squared $M^2(\sigma,\tau)$ of our approximant to Riemann Zeta function. With this Lagrangian we then proceed to apply the Least Action Principle of classical physics to get a simple and analytic action function A in the complex plane $(\sigma,\tau)$, that once plotted allowed us to find those non-trivial roots with good accuracy in a totally original way, and without resorting to algorithms or looking for changes of signs of the zeta function along the critical line. In Sub-section 5.1 we also present useful parametric plots of our Action function, in the $\sigma$-domain *(0.1, 0.9)*, using up to ten different $\tau$-coordinates of equal number of non-trivial roots as parameters. All these parametric plots significantly point to the expected common minimum at the expected abscissa value $\sigma_0=0.5$ independently of the location of the roots along the critical line chosen by us, thus confirming the validity of our application of Hamilton Principle to our approximant to Riemann Zeta. We are also able to apply this analytical procedure in a more general context using a test function based on $M^{1/\sigma}(\sigma,\tau)$, instead of the initially used modulus squared. The results we get imply necessarily that $\sigma=0.5$ is the abscissa intercept of the critical line with the real axis, as stated by Riemann Hypothesis. Our analysis also allows us in Sub-section 5.2, deriving a new function that satisfies an equation analogous to Riemann Functional Equation. Finally, in Section 6 we obtain an approximant to the Chebyshev function of second kind, and then replace it in a well-known sharp inequality estimate supported by the Riemann Hypothesis [3, 4-7] were it valid; our approximant to Chebyshev's function of second kind is also found to satisfy that sharp estimate.

## 2. A new Prime Numbers Generating Function

It is well-known that prime numbers have been studied since Euclid of Cyrene (350 BC), who showed us that there are infinite prime numbers, and how to obtain them, indeed a problem of utmost relevance in Number Theory. Some of the formulas derived thus far to obtain prime numbers are just versions of the famous *Sieve of Erathostenes* of Alexandría (ca. 280 BC) [2-4, 11]. In this section we introduce our Prime Numbers Generating Function that might yet be considered to be another version of the classical sieve, but that this time shall be found instead to be analytic, and not depending upon prime numbers tables or algorithms can be applied in all circumstances. For instance, our primes generating function allows you to explicitly calculate discreet derivatives and integrals (an example is shown in Fig. 3, Sub-section 2.3). To the effect we begin introducing below a set of required auxiliary functions: particularly, a prime numbers Discriminating Function, and a prime numbers Counting Function.

**2.1 A Prime Numbers discriminating function**

We introduce here a function that should distinguish whether a given integer is or not prime. The definition of such Prime Numbers Discriminating function, that we denote $\Lambda$ below, is therefore straightforward:

"$\Lambda(u) = 1$ if $u$ is prime, while $\Lambda(u)=0$ if $u$ is equal to 0, 1, or to any composite integer".

To obtain the expression of this new prime discriminating function we begin defining the following three simple functions of a positive real number $u$ in terms of the well-known *floor*, or *integer part*, function $\lfloor u \rfloor$:

$$\Delta(u) = u - \lfloor u \rfloor, \tag{2}$$

$$h_1(u) = 1 + \left\lfloor \frac{u}{2} \right\rfloor. \tag{3}$$

$$h_2(u,m) = 1 + \left\lfloor \frac{(u+2m-1)}{2(2m-1)} \right\rfloor \tag{4}$$

With these functions we now construct the following three auxiliary functions of $u$, and the integer's $m$ and $n$, in terms of the function $\eta(u,u_0) = sign(u-u_0)$:

$$\Omega_0(u) = [\eta(|u|,0)]^2 [\eta(|u|,1)]^2 [1 - \eta(\Delta(u),0)] \tag{5}$$

$$\Omega_1(u,m) = [\eta(|u|,2m)]^2, \tag{6}$$

$$\Omega_2(u,m,n) = \{\eta[|u|,(2m-1)(2n-1)]\}^2, \tag{7}$$

where $m > 2, n > 2$. With these three functions we may now define our prime numbers discriminating function $\Lambda$ as the double-product:

$$\Lambda(u) \equiv \Omega_0(u) \prod_{m=2}^{h_1(u)} \prod_{n=2}^{h_2(u,m)} [\Omega_1(u,m)\Omega_2(u,m,n)], \tag{8}$$

a function that appears plotted, in the integer domain *[0, 101]*, in figure 1 below:

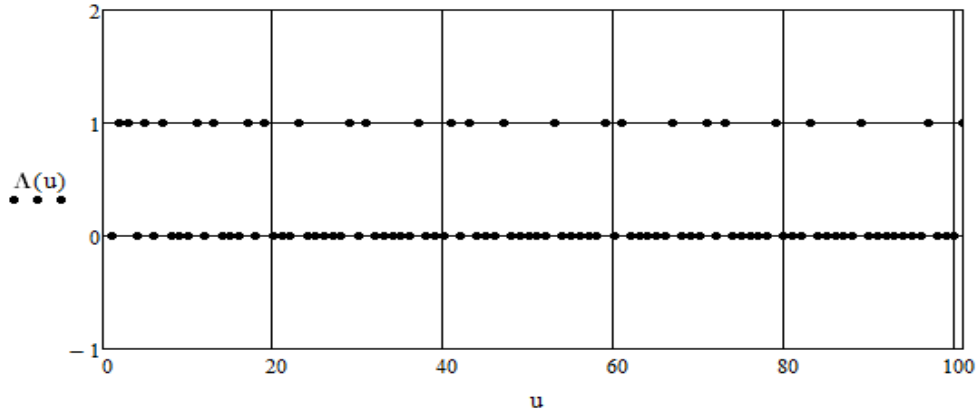

**Fig 1** Plot example of the Prime Numbers Discriminating function $\Lambda$ in the integer domain *[1,101]*,: the 26 primes (dots) appear along the straight line $\Lambda(u)=1$, the remaining integers in that domain appear on the abscissa axis.

### 2.3 A new Prime Numbers Generating function

The distribution of prime numbers among the set of integers numbers is frequently assumed to be random and has attracted mathematician's attention for centuries. Several polynomial relations to generate such numbers have been concocted, and the first that must be mentioned is the well-known polynomial formula $n^2 - n + 41$, of L. Euler (1772), that gives a list of prime numbers for $n \leq 40$. Different prime numbers generating functions have been proposed since Euler's days, some of them being recurrence formulae *e.g.* the 2008 Rowland's formula [12].

We may now proceed to define a new Prime Number Generating function $\Psi$ in terms of our prime numbers discriminating function $\Lambda$, already defined in Eq. (8). As such, our Prime's generator is a function written in terms of elementary functions, plus the *sign* function (if required the *sign* function itself may be here replaced by one approximant, defined in terms of elementary functions too, presented in our previous work [13]). Our prime generator can be used either to generate all prime numbers in a given range of integers, or actually it may also tell us whether any integer, or even a real number *u,* is or not prime. Our prime`s generator function $\Psi$ we define, in terms of our Prime Number Discriminating function $\Lambda$, as:

$$\Psi(u) \equiv u\Lambda(u) \tag{9}$$

whose first and second order *discrete* derivatives are, respectively:

$$\Psi_{d_1}(u) = \Psi(u+1) - \Psi(u) \tag{10}$$

$$\Psi_{d_2}(u) = \Psi(u+2) - 2\Psi(u+1) + \Psi(u) \tag{11}$$

Figure 3 shows a plot of this prime`s number generating function $\Psi$ for real numbers $u$ in the domain *[1,100]*: the dots of the plot correspond of course to the prime numbers in that domain; zeroes along the abscissae *u*-axis correspond of course to non-prime integers.

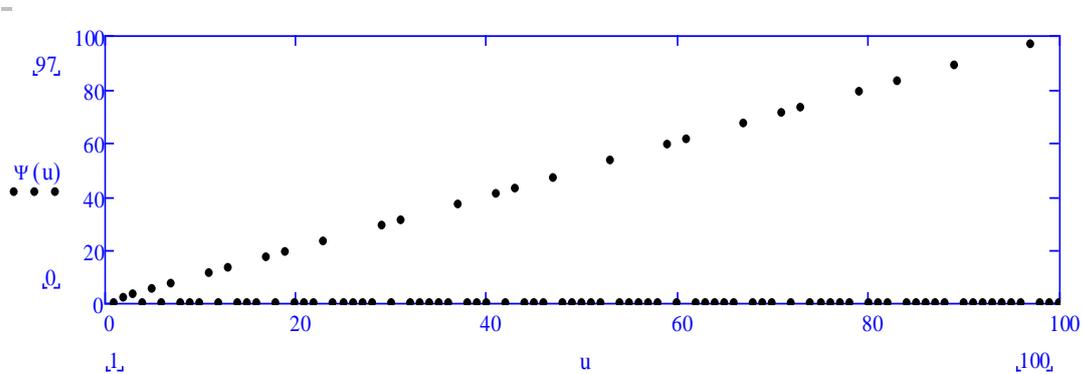

**Fig. 2** Plot of the primes numbers in the real domain (0, 100) given by our prime number generating function $\Psi$: as expected there are 25 dots lying on the inclined line of primes.

In Fig. 3 we show a plot of the first discrete derivative of our analytic prime number generating function in the domain *[1,100]*

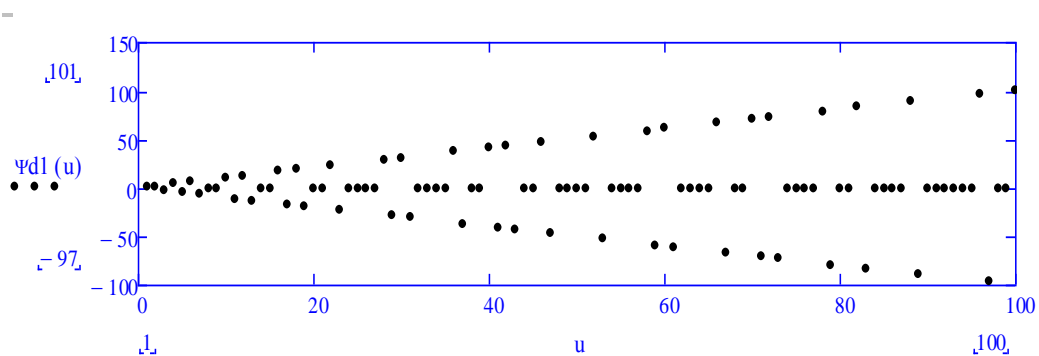

**Fig. 3** Plot of the first discrete derivative of our prime numbers generating function $\Psi$.

To emphasize the reliability of our prime number generating function $\psi$ we present Table 1 in Appendix 1, a table that shows the expected results of applying this function to 64 randomly generated integers $u_n$, that happened to fall in the integer domain *[13, 2591491]*. To the effect Table I is divided in four sets of 16 integers each, and shows the expected correct result $\psi(u_n)=u_n$ if $u_n$ is prime, otherwise returning us a zero, giving us an example of a finite but completely random distribution of primes in that integer domain. Notably, and for the sake of rigour, each of the sets of 16 random integers $u_n$ were generated using a discrete differential stochastic formalism (based on two randomly generated parameters $K_i, p_i$) as duly described in detail in Appendix 1.

**2.3 Prime Number Theorem and a new Prime Numbers Counting Function**

For the last three centuries the desideratum of Number Theory has been to obtain a relation that would give us the number of primes equal or less than any integer *x*, *i.e.* a counting function *Π(x)* that will return the exact number of prime numbers not exceeding a given *x*. In 1763 K. F. Gauss, after analyzing cases of distribution of the primes up to *large integers x*, conjectured his simple Primes Count estimate *x/ln x* in such distributions, so that we may write:

$$\frac{x}{\ln x} \approx \Pi(x) \qquad (12)$$

Years later Gauss once again conjectured [2-4] that, for large integers *x*, the number of primes *Π(x)* is also approximately given by the following function defined as a *logarithm-integral*, and known as the *Li(x)* function:

$$Li(x) \equiv \int_2^x \frac{dt}{\ln t} \approx \Pi(x) \qquad (13)$$

From Eq. (13), one may therefore write:

$$\lim_{x \to \infty} \left[ \frac{\Pi(x)}{\left(\frac{x}{\ln x}\right)} \right] = 1 \qquad (14)$$

a relation known as the *Prime Number Theorem* later proved independently by C-J de la Vallée Poussain and J. Hadamard (*ca.*1900) [2,3]. The *first-order* asymptotic $Li_a$ approximation to the *Li* function thenceforth being the first Gauss conjecture mentioned above *i.e.*:

$$Li_a(x) = \frac{x}{\ln x} \qquad (15)$$

### 2.3.1. A new Prime Numbers Counting function

Our discriminating function *Λ*, presented in Sub-section 2.1, allows us to define a new Prime Numbers Counting function $C(u, u_{in})$ which should give us the number of primes in any domain $[u_{in}, u]$ starting at a known initial $u_{in}$:

$$C(u, u_{in}) \equiv \sum_{j=u_{in}}^{u} \Lambda(j) \ . \qquad (16)$$

With this new function we may readily find the number of primes in any interval, *e.g.* for the interval *(1,1000)* we get the correct result *C(1000,2)=168*. In Fig. 4 we have plotted our *exact* prime`s number counting function *C* (the stair-case line) in the real domain *u ∈[1,10²]*. For comparison, the figure also shows the number of primes given by the first order asymptotic Gauss prime numbers counting function $Li_a(x) = x/ln(x)$ (the dotted line). As a simple comparison, the relative errors of Gauss approximate counting function and of our Prime Numbers Counting function for *x =1000* are, respectively:

$$\frac{\frac{1000}{\ln 1000}-168}{168} = -0.1383, \quad \text{and} \quad \frac{C(1000,2)-168}{168} = 0$$

Note that in Fig. 4 we have also plotted Gauss function $Li_G$. This figure illustrates that both Gauss $Li_G$ function (dotted curve) and its first-order asymptotic expansion $Li_a$ are just approximations to the exact number of primes not exceeding a given prime value, as given by our exact prime number counting function (the staircase line): ours being therefore an optimum prime counting function.

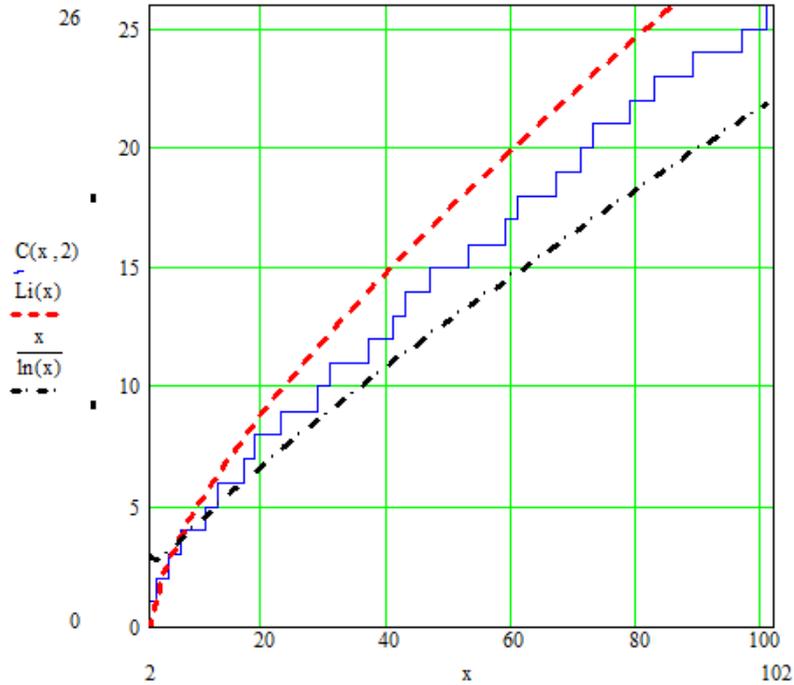

**Fig. 4** Plot (*stepped line*) of the prime number counting function *C(x,2)* in the integer domain *[1,100]*. It may be seen that for *u=100* our counting function gives the exact result *C(100,2)=25*. For comparison, also plotted appears the Gauss function *Li* (*dashed curve*) and the first order asymptotic approximation *x/ln(x)* (*dash--dotted curve*). One finds that *Li(100)=21.27* while its asymptotic approximation gives $Li_a(100) \equiv x/\ln(x) = 29.08$.

Using a modest personal computer (64 bits) and our Prime Counting Function one may readily find that the number of primes in the domain *(2, $10^5$)* is *C($10^5$,2)=9592*, which is the correct and exact result. For another comparison between those three functions plotted in Fig. 4: Gauss integral approximant gives *Li($10^5$)= 9629.62* (with *relative error ~$10^{-3}$*) for the same domain, showing it to be of acceptable accuracy. The first-order asymptotic approximation *Eq.(15),* gives $Li_a$ *($10^5$)= 8685.89*, being thus a poorer approximant to our prime counting function *C(x,2)* and of course to the ideal *π(x,2)*.

In 1901 Von Koch [14] proved that for the Riemann Hypothesis to hold then the functions difference *π(x)−Li(x)* must be of order $O(\sqrt{x}\ln(x))$. A related and important prime number theorem proved a few decades ago (1976) by Schoenfield [6**,** see also Section 6 below] includes the following corollary on a *sharp estimate* of the exact value of that same error difference:

"If Riemann hypothesis holds, then $|\pi(x) - Li(x)| < \sqrt{x} \ln x/(8\pi)$ for $x \geq 2657$"  (17a)

In figure 5 we have plotted this sharp inequality estimate of Schoenfield's corollary, in the domain *[2500, 3000]* but now replacing the prime number function $\Pi(x,2)$ with our prime counting function. The figure shows that our counting function $C(x,2)$ does satisfy this corollary related to the validity of Riemann Hypothesis. Yet, recall that our counting function has the advantage of being defined in terms of elementary functions and very easy to calculate, which makes it truly useful for its applications in relations involving prime numbers functions, as it happens in this case of the sharp estimate in Eq. (17a) above, in which we have inserted our counting functionto replace the function $\Pi(x,2)$.

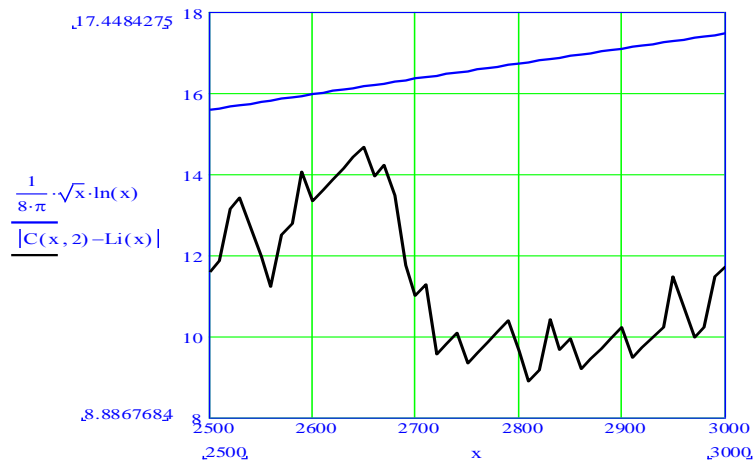

**Fig. 5** Plot of the Shoenfield's corollary sharp inequality estimate in Eq (17 a) (with ($\sqrt{x} \ln x$)/($8\pi$) plotted in blue) but replacing $\pi(x,2)$ by our Counting function $C(x,2)$, which is thenceforth shown to satisfy the inequality and proving that the Riemann Hypothesis must holds.

Moreover, just a few years ago Trudgian [15], after a thorough revision of Shoenfield's theorem, reconsidered once again the difference $\pi(x)–Li(x)$ and derived his improved version of the sharp inequality estimate for this difference that if valid proves that Riemann Hypothesis holds:

$$|\pi(x,2) - Li(x)| < \frac{0.2795\, x}{(ln(x))^{3/4}} \cdot exp\left(-\sqrt{\frac{ln(x)}{6.455}}\right). \qquad (17\ b)$$

This new inequality appears plotted in figure 6, but once again we have replaced the ideal function $\pi(x,2)$ by our analytic counting function $C(x,2)$. It may be seen that our function $C(x,2)$ does satisfy the inequality in (17 b).

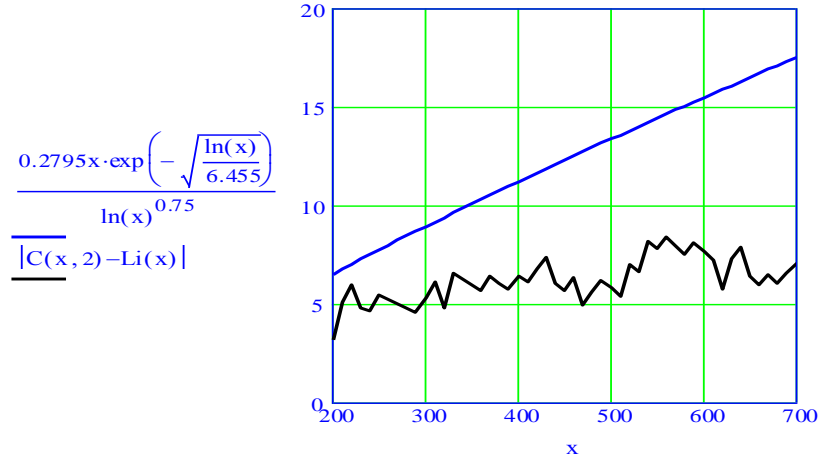

**Fig. 6** Plot of the Trudgian's corollary sharp inequality estimate, in Eq (17 b), but replacing function π(x,2) with our counting function *C(x,2),* which is thenceforth shown to satisfy the inequality, and proving that the Riemann Hypothesis must hold. (*C(2,x)-Li(x)*) appears plotted as the black line)

## 3. Euler`s Prime Numbers formula: an approximation

Let us begin this section recalling L. Euler's (*ca.*1734) remarkable real function ξ(σ), defined [2-4] for arguments σ∈ℝ, as a Dirichlet's type of convergent series:

$$\xi(\sigma) = \sum_{n=1}^{\infty}\left(\frac{1}{n^{\sigma}}\right), \tag{18}$$

which for the integers σ =2, 4 and 6 converges to the three notable values [16]:

$$\xi(2) = \frac{1}{1^2} + \frac{1}{2^2} + \frac{1}{3^2} + \cdots = \frac{\pi^2}{6}; \quad \xi(4) = \frac{1}{1^4} + \frac{1}{2^4} + \frac{1}{3^4} + \cdots = \frac{\pi^4}{90};$$

$$\xi(6) = \frac{1}{1^6} + \frac{1}{2^6} + \frac{1}{3^6} + \cdots = \frac{\pi^6}{945}$$

By invoking the Fundamental Theorem of Arithmetic, Euler later derived another formula for its real function ξ, but now remarkably *derived in terms of prime numbers p, and as an infinite product* [2], instead of the infinite sum in Eq. (18):

$$\xi(\sigma) = \prod_{p=2}^{\infty}\left(\frac{p^{\sigma}}{p^{\sigma}-1}\right), \quad \sigma > 1 \tag{19}$$

Thus after incorporating this Euler's product formula into Eq. (18) we may write:

$$\xi(\sigma) \equiv \sum_{n=1}^{\infty}\left(\frac{1}{n^{\sigma}}\right) = \prod_{p}\left(\frac{p^{\sigma}}{p^{\sigma}-1}\right), \sigma > 1 \tag{20a}$$

which we may evaluate up to a finite upper prime number limit $q_{max}(H) = P(H)$, the latter being the *H-th* prime number, to get approximated values $\xi_{eul}$ using the following Eq. 20 (b)

$$\xi_{eul}(\sigma) \equiv \sum_{n=1}^{q_{max}(H)} \left(\frac{1}{n^\sigma}\right) = \prod_{j=1}^{H} \left(\frac{P(j)^\sigma}{P(j)^\sigma - 1}\right), \quad \sigma > 1 \tag{20b}$$

Below we introduce our novel and original approximant to Euler´s product function, Eq (21), in terms of our analytic Prime Number Generating function $\Psi$ and using a finite upper limit in the product, instead of the infinite limit in Eq (20). We name this new approximant as the *self-consistent approximated Euler product* formula $\xi_{Eulap}$, and define it as the product from *2* to the $q_{max}(H)$'*th* prime number:

$$\xi_{Eulap}(\sigma, H) = \prod_{q=2}^{q_{max}(H)} \left\{ \left[\frac{\Psi(q)^\sigma + \left[\frac{\Psi(q)}{q} - 1\right]}{\Psi(q)^\sigma - 1}\right] \right\}. \tag{21}$$

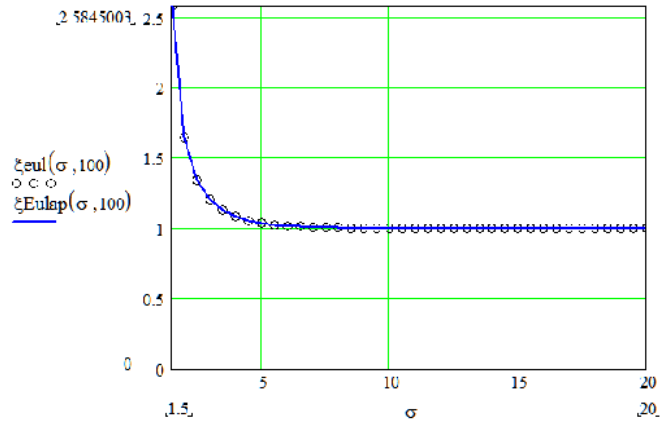

(a)

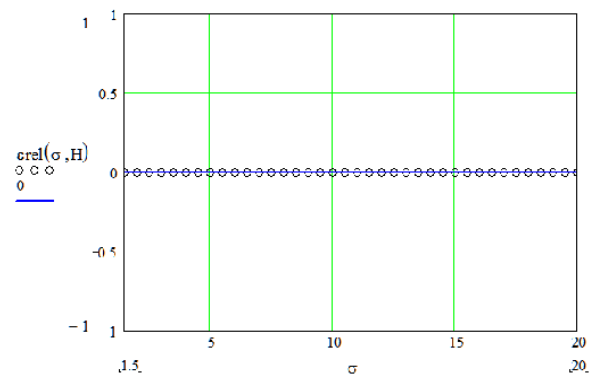

(b)

**Fig. 7 (a)** Plots of Euler's product exact function $\xi_{ex}$ and our approximant $\xi_{euapp}$ as a function of real $\sigma$. The two plots are coincident. **(b)** Plot of the relative error of our approximant to Euler's product formula (both plots for the upper bound $H=100$ of the approximant)

In this Eq. (21) the integer index $q$ in the product runs from 2 up to $P(H)$. Our Eq. (21) is a self-consistent expression that does not depend upon tables or any external algorithm, and as shown below it is a very accurate approximant. For instance the $H = 100^{th}$ prime is $q_{max}(100)=541$ and for this value our approximant *Eq.(21)* gives:

$$\xi_{Eulap}(2,100) = 1.644515221724293$$

which coincides with the exact value just obtained above using Eqs.(20). Other approximate values given by our approximant in *Eq. (21)* are:

$$\xi_{Eulap}(4,0,100) = 1.0823232233369194, \quad \xi_{Eulap}(3,0,100) = 1.202056602179509, \quad (22)$$

which compare very well with their values given by the exact Euler's product formula:

$$\xi_{eul}(4,0,100) = \frac{\pi^4}{90} = 1.082323233711138, \quad \xi_{eul}(3,0,100) \, 1.20205690315959, \quad (23)$$

respectively. A typical relative error of our approximant to Euler's product is rather small:

$$\frac{\pi^2/6 - \xi_{Eulap}(2,0,100)}{\pi^2/6} = 0.0002546 \, .$$

## 4. Approximant of the analytic continuation of Riemann function in the complex plane

A century after the works of Euler, the German mathematician Bernhard Riemann extended the application of Euler's initial formula in Eq. (18) by considering the variable *s* to be a complex number. To that effect Riemann introduced his famous Zeta Function $\zeta$, defined for any complex number $s = (\sigma, \tau)$, as the infinite sum:

$$\zeta(s) = \zeta(\sigma, \tau) = \sum_{n=1}^{\infty} \frac{1}{n^s}. \quad (24)$$

From Eqs. (20) and (24) we of course find that for real values $s = \sigma \in \mathbb{C}$ the Riemann Zeta $\zeta$ function becomes the same as Euler well-known product function $\xi_{Eul}$ over the set of prime numbers $P(j)$. As we did with Euler's Product above, we shall here search for a finite approximant to *Riemann Zeta function* considering in the present case our Prime Numbers Discriminating function $\Lambda$ defined above in Section 2.

In 1859 B. Riemann [2, 4, 11, 17-19] wrote his now famous conjecture: *The Zeta function $\zeta$ has an analytic continuation to the whole complex plane $s = \sigma \in \mathbb{C}$ except for a pole at (1,0), and it vanishes at all negative integers on the real axis (the so called trivial zeroes), and the rest of its zeroes lie on the line (1/2, $\tau$), $\tau \in \mathbb{R}$ of the complex plane.* This line is now known as the *critical line* of the continued Riemann function in that plane.

It has been found that a *restricted* analytic expression for the analytic continuation of Riemann Zeta function for all complex number *s.t. Re(s)>0* is obtained with the following complex-valued function defined as an infinite sum [2, 20, 21]:

$$\zeta(s) = \frac{1}{(1-2^{1-s})} \sum_{n=1}^{\infty} (-1)^{n-1} n^{-s}, \quad s = \sigma + i\tau, \quad Re(s) > 0 \qquad (25)$$

which is valid in the so-called *strip* of the complex plane where $\sigma$ belongs to (0,1) for any real coordinate $\tau$.

A simple and straightforward *finite* numerical expression $\zeta_{ex}$ of the infinite-sum function in Eq. (25) can be obtained by simply limiting the sum in that equation to a finite integer upper bound, say $n_{max}$=100. The real and imaginary parts of this function $\zeta_{ex}$ appears plotted in the following Fig. 8.

$$\zeta_{ex}(\sigma, \tau, n_{max}) = \frac{1}{(1-2^{1-s})} \sum_{n=1}^{n_{max}} (-1)^{n-1} n^{-s}, \quad s = \sigma + i\tau. \qquad (26)$$

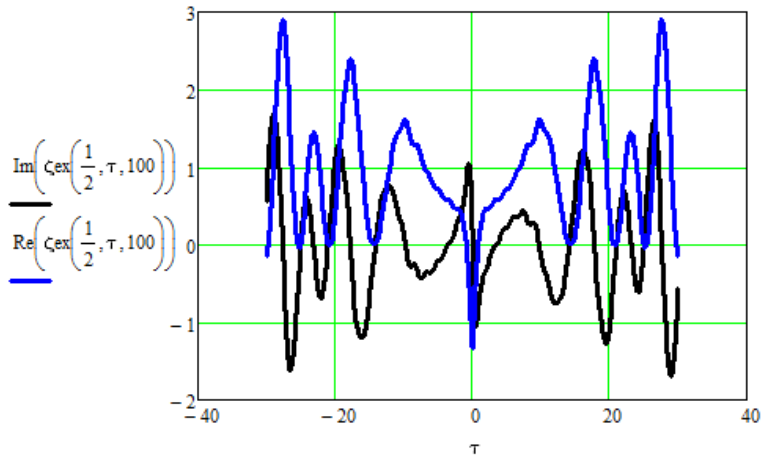

**Fig. 8** Real and imaginary part of the approximated exact Riemann Zeta function in the $\tau$-domain *(-30, 30)* along the *critical line*.

As said above we may also obtain an approximant to the Riemann Zeta function in Eq. (25) by using our Prime Discriminating Function $\Lambda$ introduced in Section 2. To the effect we simply recall that there are two classes of natural integer numbers: Primes and Composites. As defined in that section our discriminating function $\Lambda$ tells us: If an integer $n$ is a prime number it may be written as $n=n\Lambda(n)$, since then $\Lambda(n)=1$; while if $n$ is a Composite integer we get the zero value $\Lambda(n)=0$. With these two function values of Primes and Composite integers we may write the following two functions

$$\zeta_P(s) = \frac{1}{(1-2^{1-s})} \sum_{n=1}^{n_{max}} \Lambda(n)(-1)^{n\Lambda(n)-1} n^{-s}, \quad s \in \mathbb{C} \qquad (27)$$

$$\zeta_C(s) = \frac{1}{(1-2^{1-s})} \sum_{n=1}^{n_{max}} [1-\Lambda(n)](-1)^{n[1-\Lambda(n)]-1} n^{-s}, \quad s \in \mathbb{C} \qquad (28)$$

where the two complex-valued functions $\zeta_P$ and $\zeta_C$ correspond to the Primes and Composite integers, respectively. The real and imaginary parts of these two of functions appear plotted in the following two figures 9 and 10, respectively.

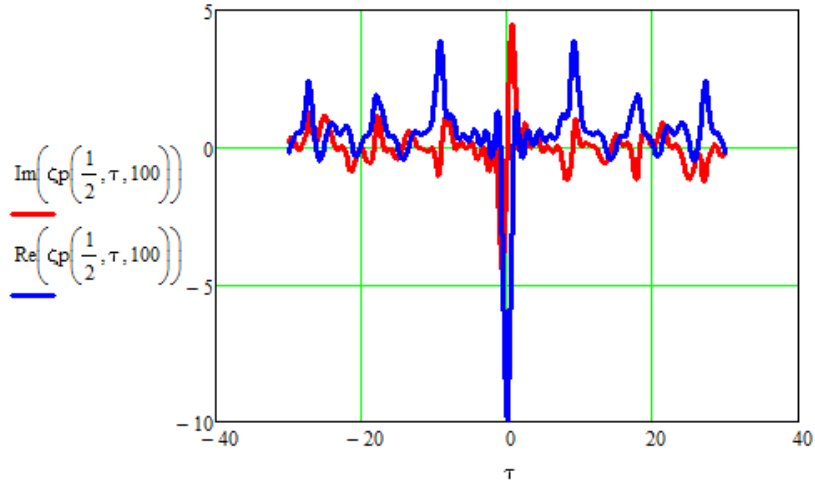

**Fig. 9** Real and imaginary parts of the function $\zeta_P$ in the $\tau$-domain *(-30, 30)* along the *critical line* (for the integer upper bound $n_{max}= 100$ in Eq. (29))

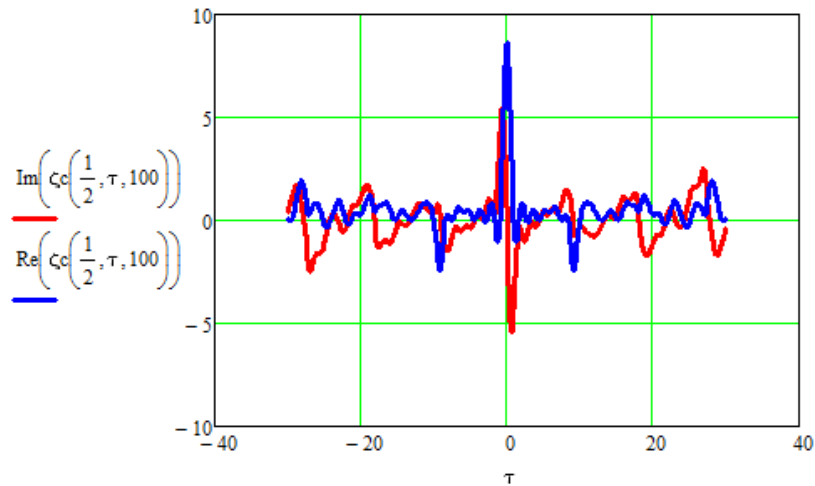

**Fig. 10** Real and imaginary part of our function $\zeta_C$ in the $\tau$-domain *(-30, 30)* along the *critical line* (for the integer upper bound $n_{max}= 100$ of the sum in Eq. (29))

We may now define a second approximant $\zeta_{app}$ to the analytic continuation of Riemann function in the complex plane, and for $\sigma>0$, as the sum of the two functions $\zeta_P$ and $\zeta_C$:

$$\zeta_{app}(s) = \zeta_P(s) + \zeta_C(s), \quad s \in \mathbb{C} \tag{29}$$

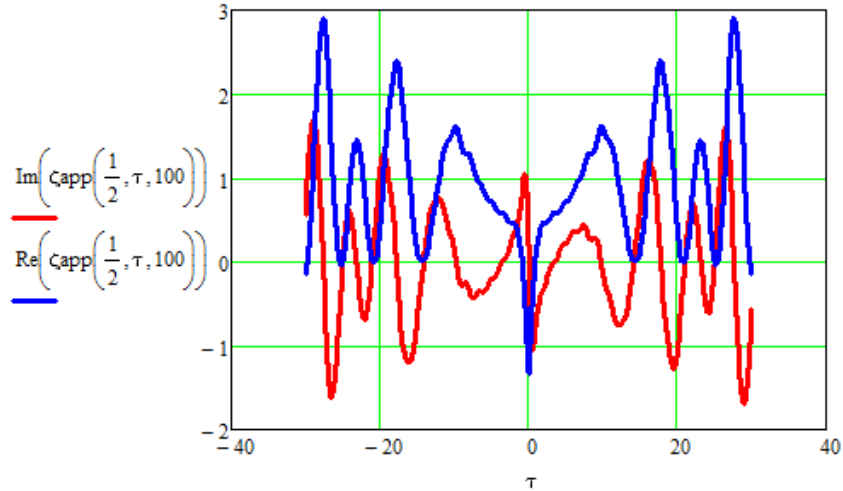

**Fig. 11** Real and imaginary part of our finite numerical approximant $\zeta_{app}$ to the Riemann Zeta function valid in the strip (0,1) with the $\tau$-domain (-30, 30), $n_{max}$=100, along the *critical line*.

The real part (blue curve) and imaginary part (red curve) of this function appear plotted in figure (11), for the upper bound $n_{max} = 100$ in Eq.(26) and in the $\tau$-domain *[-30,30]* of the Riemann critical line.

This definition of our approximant to Riemann zeta function, as a sum of the functions $\zeta_P$ and $\zeta_C$ in Eq.(29 ), is very relevant since it allows to find the approximant function values without using Prime Number tables or algorithms, but using instead our easy-to-use analytic discriminating function $\Lambda$. Moreover, it had opened to us the door to obtain an analytic formalism in which a differential equation for Riemann zeta function is derived (to be published elsewhere).

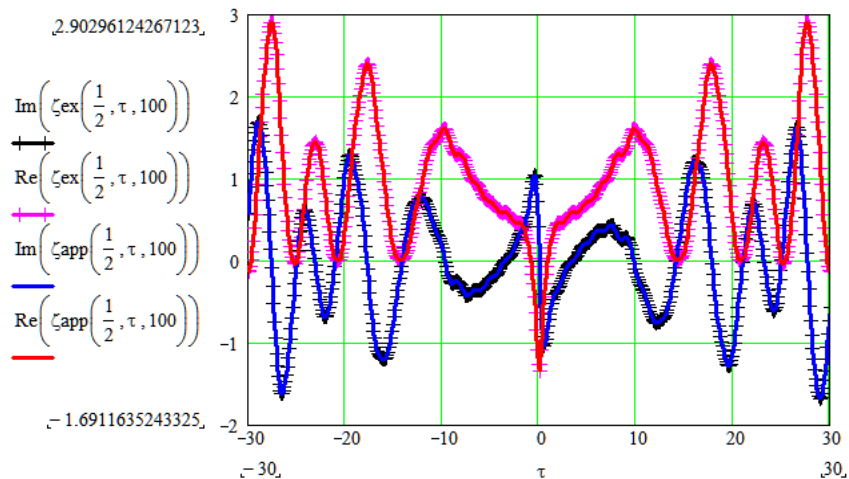

**Fig. 12** Real and Imaginary parts on the exact Riemann zeta function ($\zeta_{ex}$) and its approximant ($\zeta_{app}$) in the $\tau$-domain *[-30,30]* of the critical line, showing the coincident common zeroes of both parts.

In Fig. 12 the real part (in blue) and imaginary part (in red) of this approximant function $\zeta_{PC}$ in Eq. (29) to Riemann Zeta function, appear plotted (for the integer upper bound $n_{max}= 100$ of the sum that appears in the approximant). Also plotted are the real and imaginary part of the exact Zeta function in Eq. (26). It may be seen in Fig. 12 that the curve of our approximant lies very close to the exact function curve.

In Fig. 13 (a), (b) we have plotted the Imaginary part *vs.* Real part of the exact Riemann analytic continuation function along the critical line. If these real and imaginary parts are seen as the two coordinates of a point moving on a plane, the appealing plot in the figure may be interpreted as its dynamical trajectory travelling to the (0, 0) point, where the root is located.

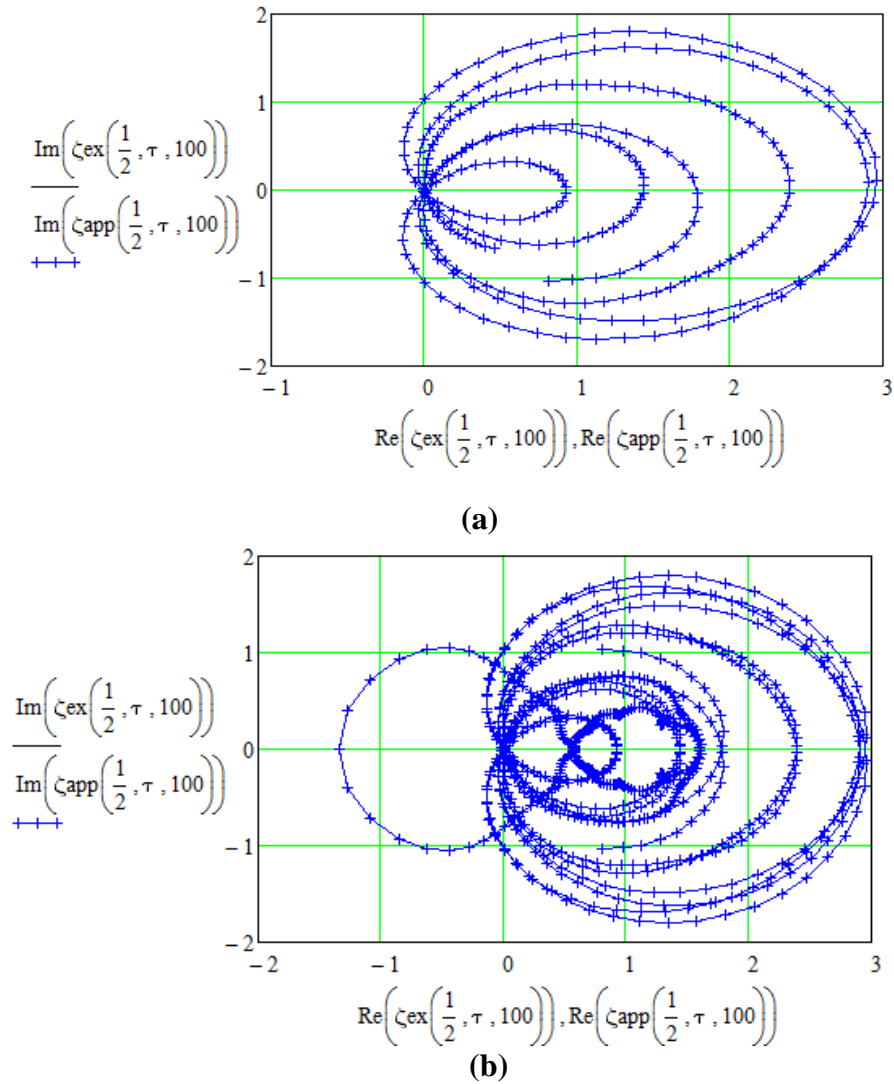

(a)

(b)

**Fig. 13** (a) (b) Plots of the real and imaginary parts of the exact Riemann function (continuous line) and of the real and imaginary parts of our approximant (crosses) to that function *vs.* ordinate ζ, respectively, along the critical line of points coordinates ($\sigma=1/2, \tau$) in the ordinate domain $\tau \in$(-0.5, 0.5).

A second and useful criterion to judge the quality of our approximant to the Zeta function is provided by plotting and comparing the modulus-squared $M^2$ of both the exact function and that of our approximant:

$$M^2[\zeta_{ex}(\sigma,\tau,n_{max})] = [Re^2(\zeta_{ex}(\sigma,\tau,n_{max}) + Im^2\zeta_{ex}(\sigma,\tau,n_{max})], \quad (30)$$

$$M^2[\zeta_{app}(\sigma,\tau,n_{max})] = [Re^2(\zeta_{app}(\sigma,\tau,n_{max}) + Im^2\zeta_{app}(\sigma,\tau,n_{max})].. \quad (31)$$

In the following figure 14 (a), (b) we have plotted these two moduli in the domain $\tau\in(-40, 40)$ of the critical line. The plotted minimae correspond of course to non-trivial roots of the analytic continuation of the Riemann function.

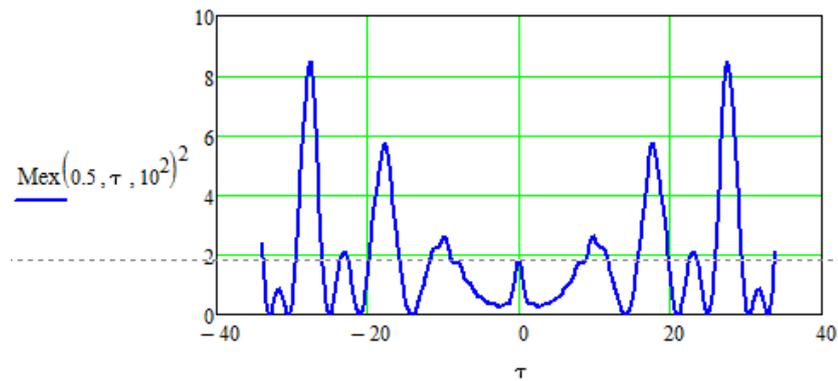

**(a)**

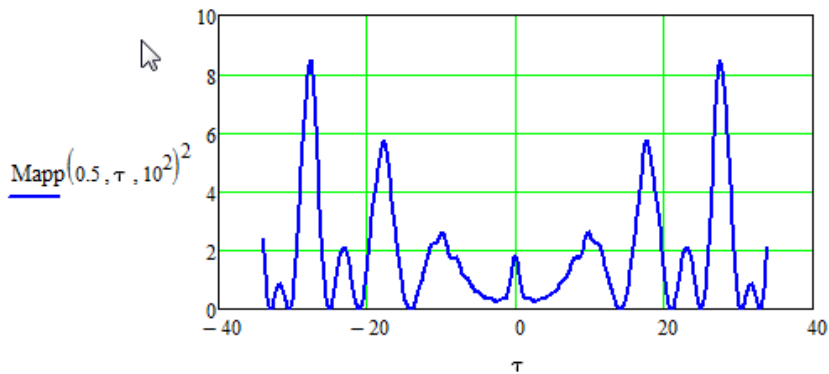

**(b)**

**Fig. 14** Comparison of the modulus-squared of the exact analytic continuation of Riemann Zeta function and of our approximant in the domain (-40, 40) of the critical line: (a) for the exact function (b) for our approximant function

In Fig. 15 we have plotted instead the reciprocal of the squared-modulus of our approximant to Riemann zeta function.

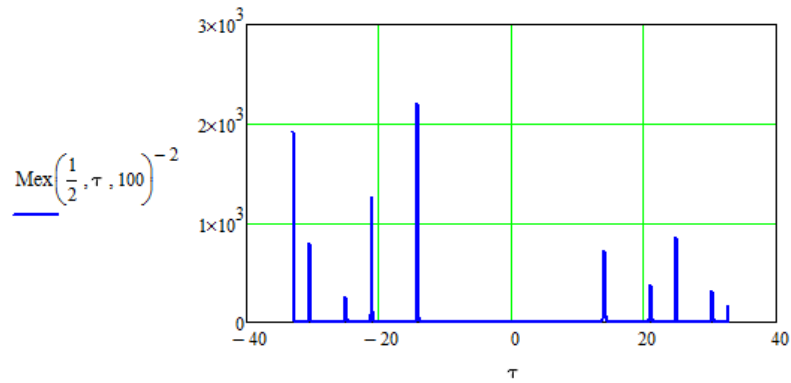

**Fig. 15** Plot of the reciprocal of the modulus-squared of the exact analytic continuation of Riemann Zeta function in the domain (-40, 40) of the critical line.

The idea of the plotting in the figure the reciprocal of the squared-modulus of our approximant is to present in a very illustrative way the accuracy of our approximant: the greater the accuracy of our approximant to a zero modulus of the Zeta function the higher, and narrower, the peak in the plot. This idea arose from the similarities of the plot in figure 15 with the familiar images one obtains when measuring the line spectra of atoms using a spectrometer in a laboratory *e.g* the well-known lines of the spectrum of light emitted by sodium atoms.

## 5. Application of the Principle of Least Action (Hamilton Principle) to the analytic continuation of Riemann Zeta function and to Riemann Hypothesis

*The Riemann Hypothesis is a precise statement, and in one sense what it means is clear, but what it's connected with, what it implies, where it comes from, can be very unobvious. (Martin Huxley)*

In this section we exploit a basic tool of classical mechanics, Hamilton Principle, [22, 23] to locate the non-trivial roots of Riemann zeta function in the complex plane, and what is more important to show that given a set of those non-trivial roots they must lie on the so-called critical line, line parallel to the τ-axis that intersects the real axis of the complex plane at the point of abscissa $\sigma_0=1/2$ as prescribed by the Riemann Hypothesis. We also apply the least action principle to derive a complex function that satisfies an equation analogous to the functional equation presented by Riemann in his hypothesis.

### 5.1 Application of the Principle of Least Action to find the non-trivial roots of Riemann zeta function

In order to apply the Least Action Principle we consider, for the sake of simplicity, the Lagrangian of an ideal harmonic oscillator, that has unit mass, unit elastic constant and unit amplitude *i.e.* an oscillator whose position is represented by *x(t)*= cos (*t*). This simple system is to be perturbed by a test function that we chose to define in terms of the modulus-squared of the Riemann's zeta function times the product *sin(t) cos(t)*, that is: $M_{ex}^2(\sigma, \tau, n_{max})\, cos(t)\, sin(t)$. In this way the variational perturbation becomes zero at the extremes of the motion time interval $(0, 2\pi)$, as it is expected to be. The variational test function $X$ of the mechanical system we are considering is therefore:

$$X(t, \sigma, \tau, n_{max}, t_i, t_f) = \cos t \cdot [1 + M_{ex}^2(\sigma, \tau, n_{max}, t_i, t_f) \cdot \sin t], \quad (32)$$

where we shall set $n_{max}=100$, $t_i=0$ and $t_f= 2\pi$ (in Eq. (34) below). The Lagrangian $L = T-U$ we shall use is therefore:

$$L(t, \sigma, \tau, n_{max}, t_i, t_f) = \frac{1}{2}\left[\frac{d}{dt}X(t, \sigma, \tau, n_{max}, t_i, t_f)\right]^2 - \frac{1}{2}X^2(t, \sigma, \tau, n_{max}, t_i, t_f)^2, \quad (33)$$

which now allows us to write the pertinent Action Integral:

$$A(\sigma, \tau, n_{max}, t_i, t_f) = \int_{t_i}^{t_f} L(t, \sigma, \tau, n_{max}, t_i, t_f)\, dt \quad (34)$$

Having defined this Action Integral, for the application of the variational calculus, we may proceed to explicitly calculate the values of the variables $\sigma$ and $\tau$, for the Action to reach a minimum, as expected from the application of Hamilton Principle. This calculation may be simplified by keeping constant the upper bound value to $n_{max}=100$ when evaluating the modulus squared of the Zeta function in Eq. (34), and by introducing a simpler notation for that modulus: $M(\sigma, \tau) = M_{ex}(\sigma, \tau)$

Analytical integration of the Action in Eq. (34) returns us the following useful and *exact analytic* expression:

$$A = \frac{3}{8}\pi M^4(\sigma, \tau) \quad (35)$$

Taking the logarithm of this expression we get a relation that shall immediately prove to be very useful:

$$log[A(\sigma, \tau)] = log\left(\frac{3}{8}\pi\right) + 4\, log[M(\sigma, \tau)]. \quad (36)$$

In effect, plotting the *log (A)* function in Eq. (36) vs. *log( M)* we obtain the linear log-log plot in Fig. 16 of the Action against the Modulus of the Zeta function in the complex plane:

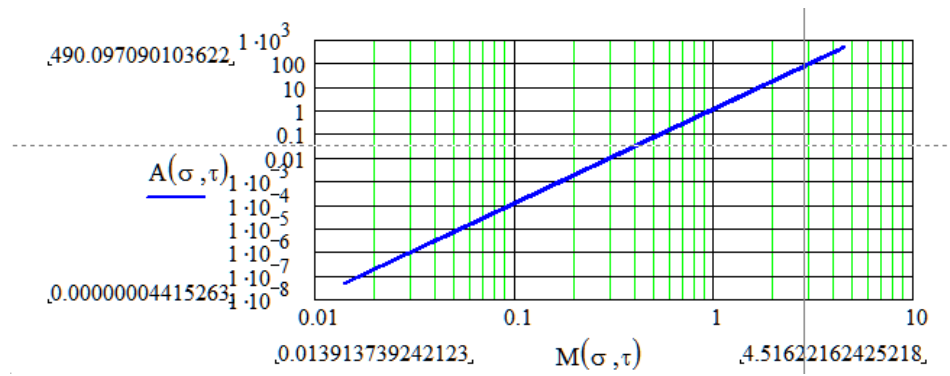

**Fig. 16** Log-log linear plot of the Action *(σ, τ)* versus the Modulus $M(\sigma, \tau)$ of the Zeta function

For further plotting of the Action function, in which this function would display for us its expected minimum values, it is convenient to define the following simple set of two variables:

$$\eta(\sigma, \tau) = \tau + \sigma, \qquad \omega(\sigma, \tau) = \tau - \sigma \qquad (37)$$

Thus in the following Fig. 17 (a) we present a semi-log plot of the analytically obtained Action, in Eq (35), against the variable $\omega = \tau - \sigma$, in the domain *(10, 45)*, to see whether we get a set of minimae. As shown in the figure the resulting plot consists of a bundle, or family, of 9 curves all pointing down to minimae values located on the abscissa axis. Each bundle of curves observed in the figure corresponds to 8 discretely varying values of the variable σ, in steps *Δσ=0.1* along the domain *(0.1,0.9)*, and with the variable τ varying from 13.0 to 44.0 in steps of 0.1.

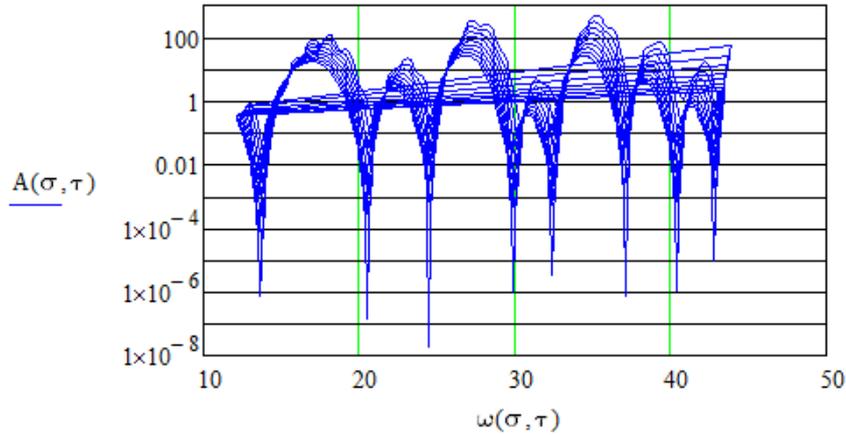

**Fig. 17 (a)** Plot of the log of the Action versus the variable $\omega(\sigma, \tau) = \tau - \sigma$ sum of the coordinates of the complex number *s=(σ, τ)*. A family of curves that points to eight very small minimae ($<10^{-6}$). Each family of curves in the figure correspond to discretely varying values of the variable σ in steps *Δσ=0.1* in the domain *(0.1,0.9)* and τ in steps of 0.1 in the range 13.0 to 44. Note the small values ($<10^{-5}$) of the action A at its minimae

Analogously, in the following figure Fig. 17 (b) we plot the analytically obtained Action in Eq (36) this time against the variable $\eta = \sigma + \tau$, in the domain *(10, 45)*. As shown in the figure the plot once again consists of a family of curves that as expected point down to minimae values located on the abscissa axis.

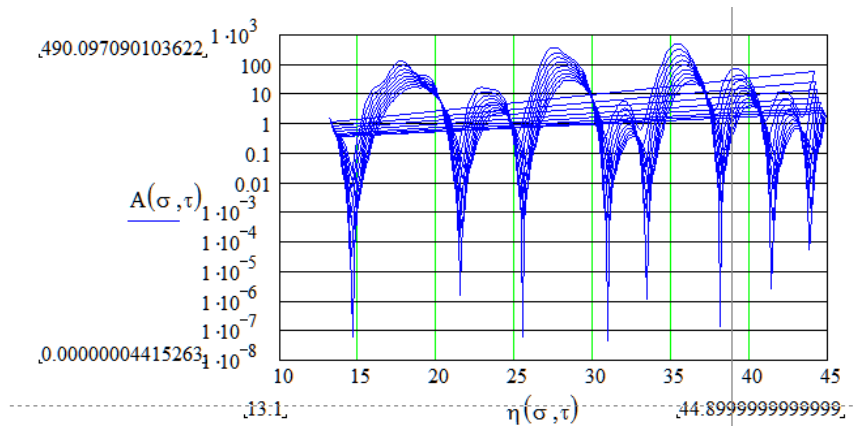

**Fig. 17 (b)** Plotting the log of the Action versus the variable $\rho(\sigma,\tau) = \sigma+\tau$ of the coordinates of the complex number $s=(\sigma,\tau)$. A bundle of curves that points to eight very small minimae (<$10^{-6}$). Each manifold of curves in the figure correspond to discretely varying values of the variable $\sigma$ in steps $\Delta\sigma=0.1$ in the domain *(0.1, 0.9)*, and $\tau$ in steps of 0.1 in the range 31.0 to 31.9

To obtain the actual coordinates *($\sigma,\tau$)* of a minimum of the Action A from Figs. (17) we may simply zoom into a small neighbourhood of any minimum, say about the minimum located by $\omega= 32$ in Fig. 17 (a). The result of this zooming is shown in Fig. 18 (a) and gives us a minimum at $\omega= 32.406$. This zooming procedure is applied again to Fig. 17 (b) and the result $\eta=33.40$ appears in Fig. 18 (b).

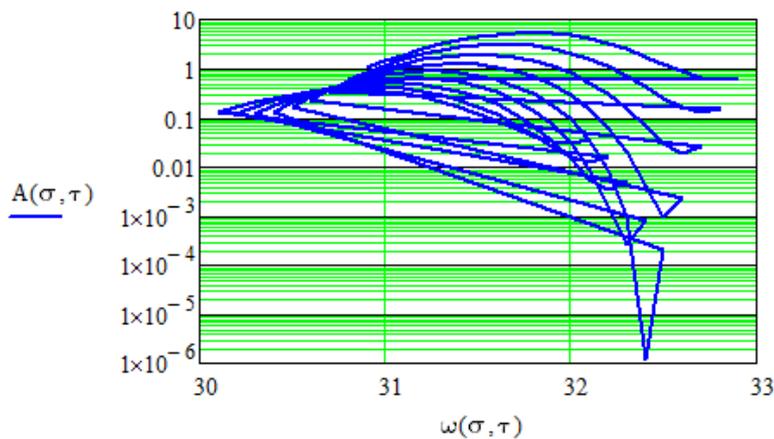

**Fig. 18 (a)** Plot of the log of the action versus the variable $\omega(\sigma,\tau) = \tau - \sigma$, sum of the coordinates of the complex number $s=(\sigma,\tau)$ in the domain *(30,33)*. The minimum occurs at $\omega=32.406$.

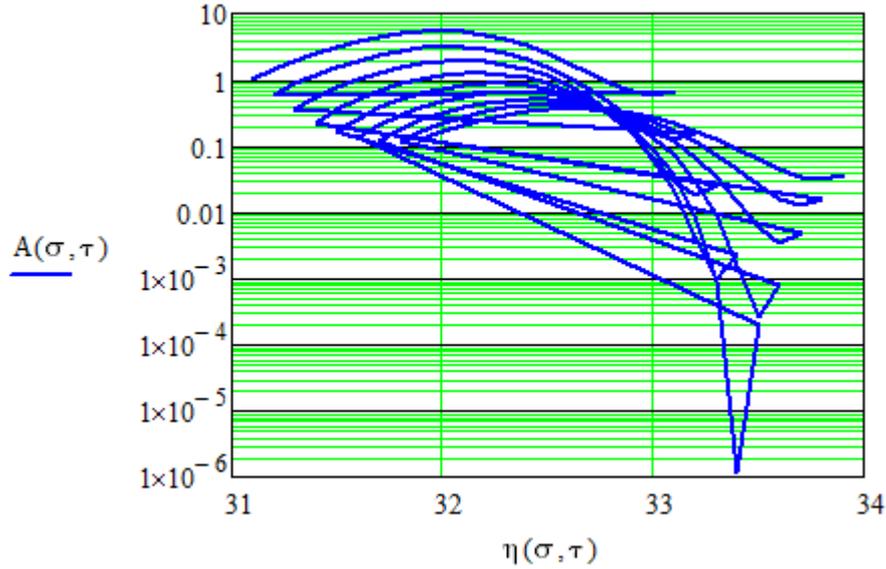

**Fig. 18 (b)** Plot of the log of the action versus the variable $\eta(o,\tau) = \sigma + \tau$, difference of the coordinates of the complex number $s=(\sigma,\tau)$ in the domain *(30,33)*. The minimum occurs at $\eta=33.4$.

From the minimae of the analytical Action obtained in figures, 18 (a) and 18 (b), it is now possible to obtain the coordinates *($\sigma,\tau$)* of the actual minimum of the Zeta function in the pertinent region of the complex plane. In effect solving the system of equations for the two minimae given by the zooming:

$$\tau + \sigma = 33.40, \quad \tau - \sigma = 32.406$$

we get the coordinates: $\sigma= 0.497$, $\tau=32.903$. The known published coordinate values of the pertinent non-trivial zero of the Zeta function are in the present case are: $\sigma=0.500$ and $\tau=32.903$. The relative errors for the found coordinates of this particular non-trivial root are 0.0060 and 0.00097, respectively. These small relative errors show the good accuracy of our method based on Hamilton Principle of physics. In fact we have successfully applied it to not less than 20 non-trivial roots of Riemann Zeta function. For instance from well-known published of non-trivial roots of Zeta [25] we know that there is one of such roots located on the critical line at the coordinates $\sigma= 0.5$, $\tau = 998.827547137$. Using our procedure described above we obtained the rather accurate result $\sigma= 0.495$, $\tau=998.83$ (see Appendix 3).

Below, in figure 19(a) a parametric plot, we plot the exact Action function we obtained in Eq. (35) as a function of the real coordinate $\sigma$ along the abscissa axis, yet using as parameter the coordinate $\tau$ of the following five different non-trivial roots of Riemann Zeta function on the critical line:

$$\rho_I \equiv \tau_I = \{14.13, 21.02, 25.01, 30.43, 49.55\},$$

and evaluating the Action using $n_{max}=100$ as the upper bound in our approximant to Riemann zeta function (Eq. 29). It may be seen that the five parametric plots show the expected single minimum of the Action to occur at the coordinate $\sigma_0=0.5$, validating in a simple but truly significant way our application of the Least Action Principle to Riemann Zeta function.

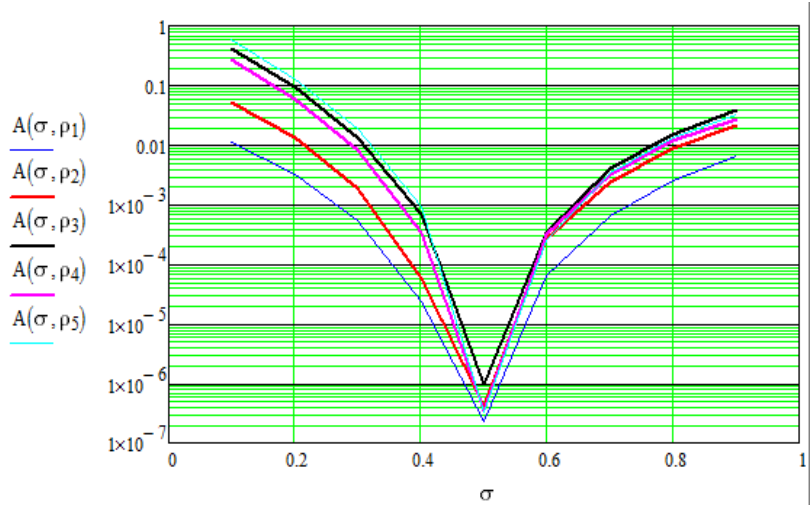

**Fig. 19 (a)** Parametric plots of the Action as a function of real coordinate $\sigma$ in the domain *(0, 0.9)*, using as parameter five different nontrivial roots $\rho_i \equiv \tau_i = \{14.13, 21.02, 25.01, 30.43, 32.94\}$ of the Zeta function along the critical line. The five different curves show the same expected minimum at $\sigma_0=0.5$.

In figure 19(b) we have again plotted the Action function, obtained in Eq. (35), as a function of the real coordinate $\sigma$, but this time using as a parameter the coordinates $\tau$ of five different but higher non-trivial roots of Riemann Zeta function along the critical line namely: $\rho_i \equiv \tau_i = 996.205, 997.511, 998.827, 999.792, 1001.349$. As expected once again we get the important result: a single minimum at the abscissa $\sigma=0.5$.

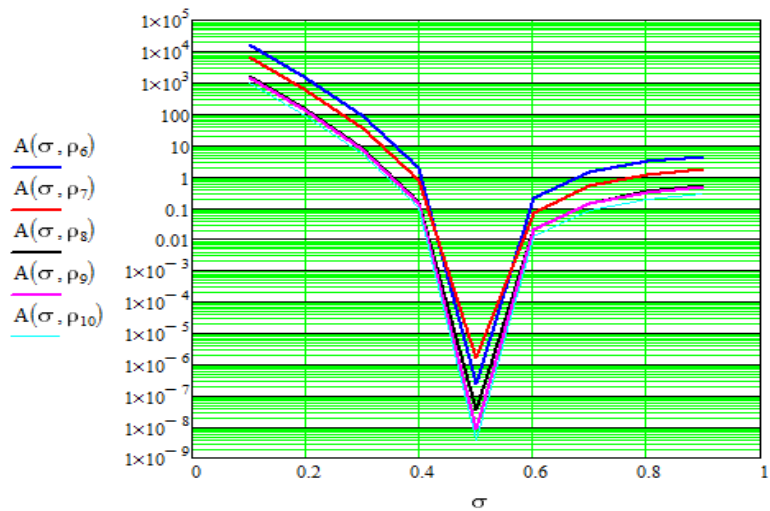

**Fig. 19 (b)** Parametric plots of the Action as a function of real coordinate $\sigma$ in the domain *(0, 0.9)*, using as parameter five different nontrivial roots $\rho_i \equiv \tau_i = \{996.205, 997.511, 998.827, 999.792, 1001.349\}$ of the Zeta function along the critical line. The Action was evaluated using $n_{max} = 10^2$ as the upper bound in our approximant to the function. The five different curves show the same expected minimum at $\sigma_0=0.5$.

Figures 19(a) and 19(b) may be unavoidably seen as a strong argument in favour of that the Riemann Hypothesis holds. In effect, it cannot be a fortuitous result that ten or more parametric plots of the Action, obtained by us, point correctly to the expected abscissa value $\sigma = 0.5$ where the critical line of non-trivial roots, according to Riemann, should intersect the real axis, independently of the ten parameter values.

**5.2 Application of the Least Action Principle and the variational calculus to obtain a Riemann-type Functional equation and non-trivial roots of the Zeta function in finite domains.**

In the previous sub-section we applied the Principle of Least Action to get non-trivial roots of Riemann zeta function. To the effect we considered a variational test function written, using simple heuristic arguments, in terms of the squared modulus $M_{ex}^2$, of the Zeta function. The mechanical system we use, for the sake of simplicity, is a harmonic oscillator of unit mass and unit elastic constant. Here we want to extend that analytical mechanics treatment, but this time using a variational test function written in terms of $(M_{ex})^{1/2}$, and applying again that analytical mechanics principle to obtain the minimum of the resulting Action function. It might then happen that the minimum value of the action would occur at a value of $\sigma$ different from the value demanded by Riemann Hypotesis, i.e. different to $\sigma = 1/2$, which will tell us that that hypothesis does not hold. Thus, let the variational test function be given by Eq. (32) but with the squared modulus of Zeta function replaced by $M(\sigma, \tau, n_{max})^{1/\sigma}$. After applying the least action principle, to the corresponding Lagrangian, the new analytical Action function that is, once again, analytically obtained is given in Eq, (38) below,

$$A(\sigma, \tau, n_{max}) = \frac{3}{8}\pi\, M(\sigma, \tau, n_{max})^{2/\sigma}, \tag{38}$$

We also need to consider the *energy dispersion $E_a$* of the mechanical system, again assumed to have unit mass and unit elastic constant, with respect to its mean value 1/2:

$$E_a(\sigma, \tau, n_{max}) - \frac{1}{2} = \frac{5}{16}[M_{ex}(\sigma, \tau, n_{max})]^{2/\sigma}. \tag{39}$$

Below in Fig. 20 we have plotted the logarithm of the modulus $M_{ex}(\sigma, \tau, n_{max})$ of the exact Riemann zeta function [2, 20,21] against the logarithm of the action $A(\sigma, \tau, n_{max})$, for nine different values of $\sigma$. The interesting plot consists of a pencil of straight lines and different slopes *m*, as expected from the application of the Least Action Principle, each line corresponding to a different value of $\sigma$. In the figure we have also plotted just *two points* of coordinates $M_{ex}(\sigma_0, \tau_1, n_{max})$ and $A(\sigma_0, \tau_2, n_{max})$, for the particular values $\tau_1 = 14.134725142$, and $\tau_2 = 21.022039639$ that correspond to the particular
$\tau$-coordinates of two non-trivial roots of function zeta, namely $(0.5, \tau_1)$, $(0.5, \tau_2)$. The figure shows that the two points in question (plotted as *circles*) belong to the same line.

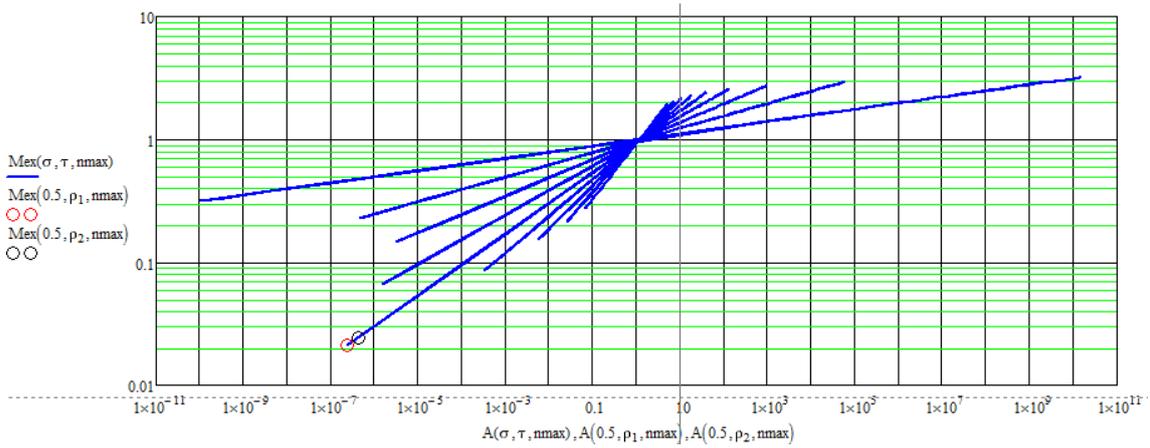

**Fig. 20** Log-log plot (pencil of blue lines) of the modulus $M_{ex}(\sigma\tau,n_{max})$ of Riemann zeta function versus the action $A(\sigma,\tau,n_{max})$ with $n_{max}=100$. The two *circles* correspond to two known roots of Zeta function at the points: $(0.5, \tau_1), (0.5, \tau_2)$.

The slope *m* of that particular line can be estimated from the figure itself just using the coordinates of the two points (circles) in Fig. 20. The result being: *m= 0.25001201719599*. On the other hand the expected slope in the log-log plot can also be readily calculated, using Eq. (39), and it is found to be: $\sigma_0/2= 0.5/2= 0.25$. Thus our results obtained applying analytical mechanics are definitely consistent.

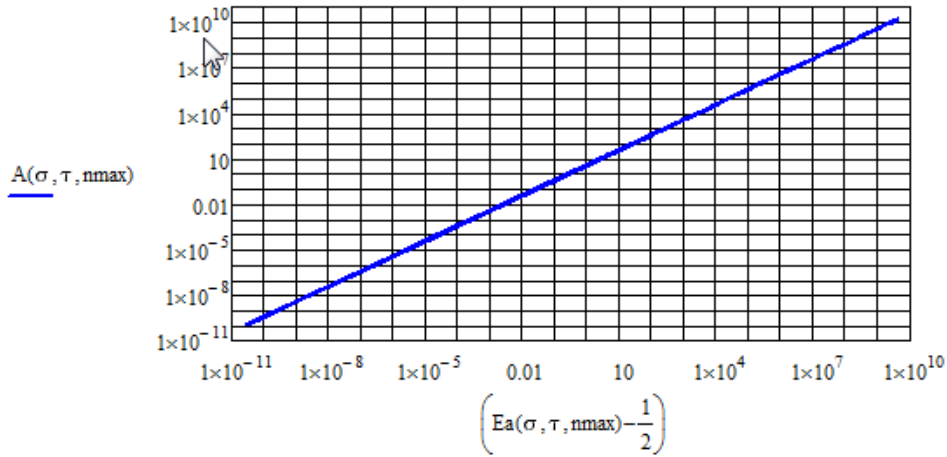

**Fig. 21** Log-log plot of the action $A(\sigma, \tau, n_{max})$ versus the energy $E_a(\sigma, \tau, n_{max}) - \frac{1}{2}$

Figure 21 is a log-log plot of the analytically calculated Action (in Eq. (39)) versus the mechanical system energy dispersion $E_a$. The plot shows a striking straight line that contains key useful information, as immediately explained below. The slope of that straight line is given by:

$$\frac{E_a(\sigma,\tau,n_{max})-\frac{1}{2}}{A(\sigma,\tau,n_{max})} = \frac{5}{16} \cdot \frac{8}{3\pi} = \frac{E_a(1-\sigma,\tau,n_{max})-\frac{1}{2}}{A(1-\sigma,\tau,n_{max})}. \tag{40}$$

Below, in Fig. (22) we plot the two ratios that appear on the left side (in *circles*) and on the right side (the *red line*), in our Eq. (40), respectively. Both are plotted as function of the relevant real abscissa coordinate $\sigma$ of the complex plane, in the domain

*(0.1,0.9)*. As seen, and strikingly, the two ratios that have been evaluated independently, take the same value for any value of $\sigma$. Note that the constant value $C=40/48\,\pi$ of the two ratios is also plotted in Fig. 22 as the small *black dots* (inside the circle).

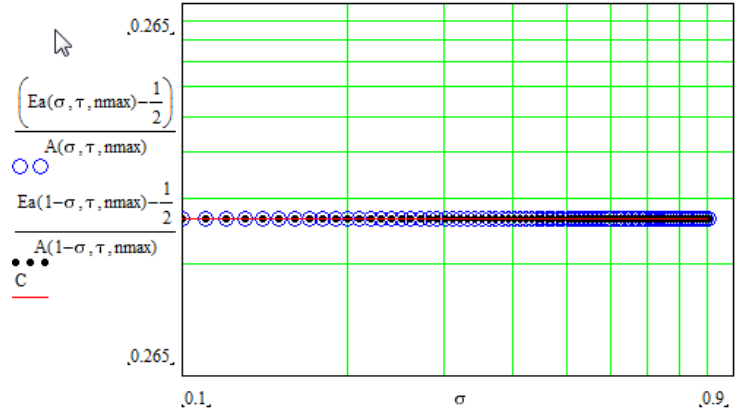

**Fig. 22** Plots of the two ratios appearing in in Eq. (41): ratio $[E_a(\sigma,\tau,n_{max}) - \frac{1}{2}]/A(\sigma,\tau,n_{max})$ (in blue circles) and ratio $[E_a(1-\sigma,\tau,n_{max}) - \frac{1}{2}]/A(1-\sigma,\tau,n_{max})$ (red line). Their common constant value $C=(5/16)(8/3\pi)$ also plotted (as black dots)

Therefore from Eq. (40), and from the plots in Fig. 22, we may now write:

$$\left[E_a(\sigma,\tau,n_{max}) - \frac{1}{2}\right] A(1-\sigma,\tau,n_{max}) = A(\sigma,\tau,n_{max}) \cdot \left[E_a(1-\sigma,\tau,n_{max}) - \frac{1}{2}\right] \quad (41)$$

This is a symmetrical relation that contains the required information for the Riemann Hypothesis to hold. In effect, this Eq (41) suggests us to introduce in our work a new function $F$ which we define below as the product of the Energy dispersion times the Action:

$$F(\sigma,\tau) \equiv \left[E_a(1-\sigma,\tau,n_{max}) - 1/2\right] \cdot A(\sigma,\tau,n_{max}). \quad (42)$$

This function happens to satisfy the symmetry property written below in Eq. (43):

$$\boldsymbol{F(\sigma) = F(1-\sigma)}, \quad (43)$$

a relation that is simply analogous to the well-known symmetrical Riemann Functional Equation [2]

The following figure 23 (a) is a semi-log plot of function $F(\sigma,\tau) = A(\sigma,\tau,n_{max}) [E_a(\sigma,\tau,n_{max}) -1/2]$, versus the variable $\tau$ along the abscissa axis. The plot shows three minimae precisely located at three known non-trivial roots of the Riemann analytic continuation function. In the present case the plot points to minimae at: 14.135, 21.035 and 25.035. The corresponding Riemann tabulated values for those roots are 14.134725142, 21.022039639 and 25.01085758 [25].

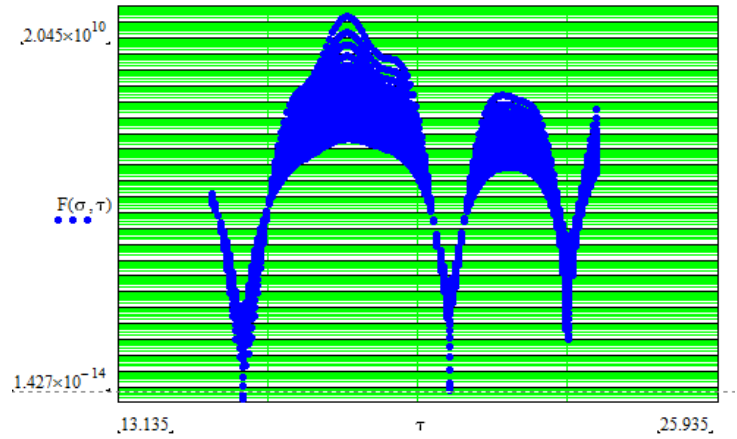

**Fig. 23(a)** Semi-log plot of the function $F(\sigma,\tau) = A(\sigma,\tau,n_{max}) [E_a(\sigma,\tau,n_{max}) - ½]$ against the imaginary coordinate $\tau$ of the complex plane showing minimae at three non-trivial roots of the Riemann zeta function in the $\tau$-domain (13.135, 25.935), namely: 14.134725142, 21.022039639, 25.01085758.

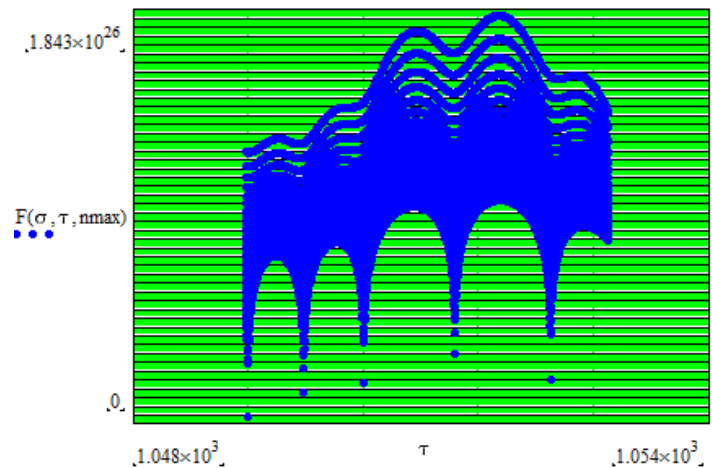

**Fig. 23(b)** Semi-log plot of the function $F(\sigma,\tau) = A(\sigma,\tau,n_{max}) [E_a(\sigma,\tau,n_{max}) - ½]$ against the imaginary coordinate $\tau$ of the complex plane showing minimae at the five non-trivial roots of the Riemann zeta function in the $\tau$-domain (1048, 1054), namely: 1047.987147490; 1048.953785195; 1049.996284257; 1051.576571843; 1053.245785158

Figure 23 (b) shows a similar plot to that in Fig. 23(a), but for a higher domain along the critical line showing minimae at five non-trivial roots of Riemann zeta function.

In figures 24(a) and 24(b), below, we have plotted instead the analytical expression of the functions $F(\sigma)$ and $F(1-\sigma)$ against the real abscissa $\sigma$ of the complex plane. The plots show the expected single minimum at $\sigma=0.5$.

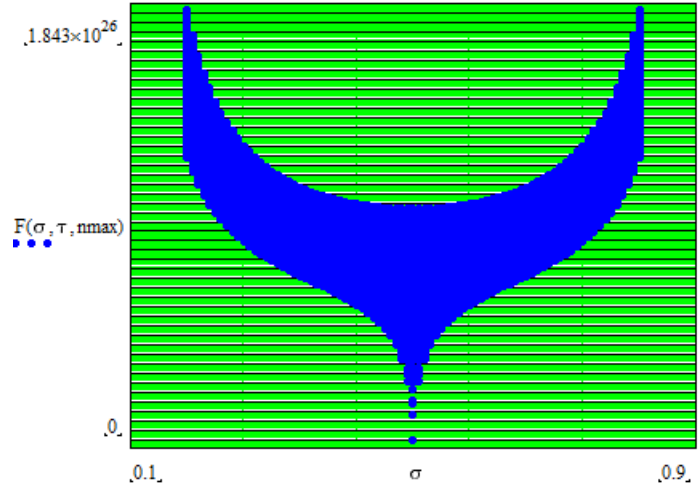

**Fig. 24 (a)** Semi-log lot of the function $F(\sigma,\tau)$ versus the real abscissa $\sigma$ of the complex plane showing the expected minimum value $F(\sigma,\tau,n_{max})=0$ at $\sigma=0.5$.

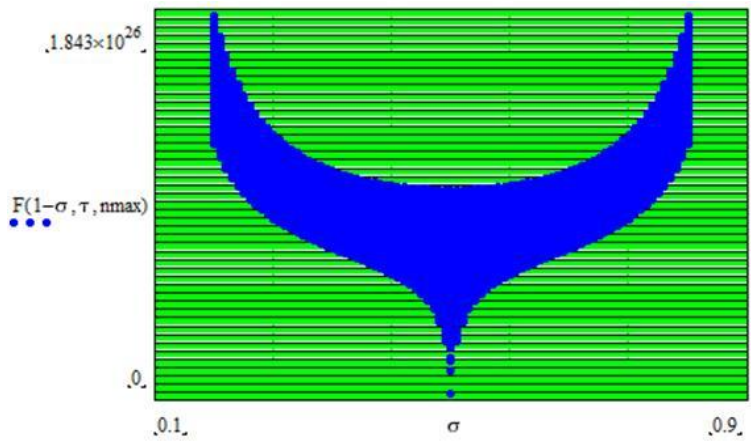

**Fig. 24 (b)** Semi-log lot of the function $F(1-\sigma,\tau)$ versus the real abscissa $\sigma$ of the complex plane showing the expected $F(1-\sigma,\tau,n_{max})=0$ at $\sigma=0.5$.

Looking for an explanation of the single minimum in both of the figures (24) we finally arrived to the expected explanation: Figures 24(a), (b) are just a direct consequence of the analogy of our function $F$ with the Functional Equation that Riemann presented in his historical mathematical academic dissertation [2, 21]. In effect: the function $F$ satisfies an analogous equation to the well-known Functional Equation of Riemann Hypothesis, and gives us the expected single minimum at $\sigma=0.5$.

## 6. Approximant to the Chebyschev Function of Second Kind

Equation (44) below is an *exact* formula for evaluating the well-known Chebyshev Function of second kind [24] in terms of the prime function $P(j)$, and the *integer part* function. The $P(j)$ function simply gives the *j-th* prime, e.g. $P(1)=2$ and $P(100)=541$.

$$\Psi_{Chex}(x) = \sum_{j=1}^{\Pi(x,2)} \left\{ floor\left[\frac{\ln(x)}{\ln[P(j)]}\right] \ln[P(j)] \right\}, \tag{44}$$

In this equation the upper limit $\Pi(x,2)$ of the sum, is the exact prime counting function i.e. the number of primes not exceeding a given integer threshold *x*. Here, instead of the function *P(j)*, we replace it with our exact Prime Numbers Generating function $\psi$ and our Prime Discriminating Function $\Lambda$, both defined in Section 2. This gives us a very useful analytical approximant function $\Psi_{Cheap}$ to Chebyshev's Function, namely:

$$\Psi_{Cheap}(x) = \sum_{q=2}^{x} \left\{ floor\left[\frac{\ln(x)}{\ln[\Psi(q)+e(1-\Lambda(q))]}\right] \ln[\psi(q) + e(1 - \Lambda(q))]\, \Lambda(q)\right\}, \quad (45)$$

where *floor* is the integer part function, and *e* is the base of the natural logarithms. Note also the finite upper limit of the sum in our approximant, in Eq (45), instead of the function $\Pi(x)$. Thus our approximant $\Psi_{Cheapp}$ is a self-consistent expression that *does not require the use of tables of prime numbers*. The good accuracy of this approximant to Chebyschev function can be appreciated considering the values, given by the two functions (44) and (45) for the *j*=100[th] prime. Up to an accuracy of $10^{-14}$ we get the same results:

$$\Psi_{Cheex}(j) = 94.0453112293574 \quad \text{and} \quad \Psi_{Cheap}(j) = 94.0453112293574$$

In figure 25 (a) and 25 (b), and for the sake of comparison, we have plotted both the exact Chebyshev function of Second kind and our approximant function to this function, in the domain *[0,200]*, respectively. It may be seen that our approximant does reproduce the plot of the exact function. This assertion is validated by our following figure 26.

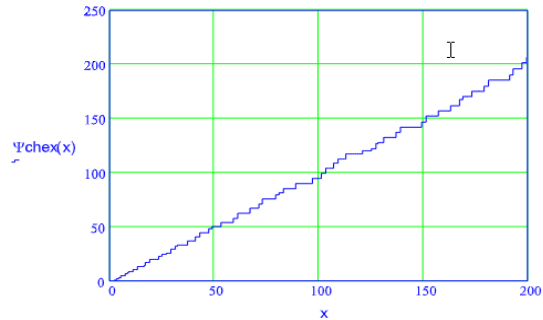

(a)

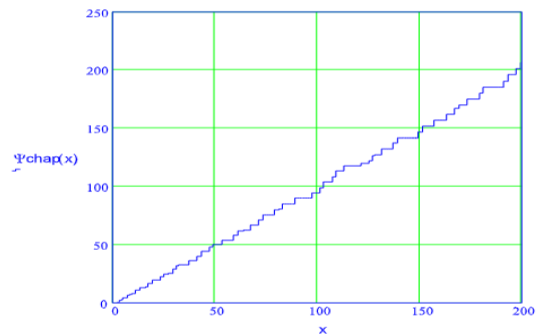

(b)

**Fig. 25** Comparison plots of the exact Chebyschev 2$^{nd}$ kind function **(a)** and of our approximant to that function **(b),** Eq. (45), in the real domain *[0,200]*.

In figure 26 we have plotted the relative error of our approximant $\Psi_{Cheapp}$ when compared with the exact Chebyshev function in the domain *[0,200]*. It may be seen that the accuracy of our approximant (Eq. (45)) which uses our Prime Number function $\Lambda$ is very high, the relative error being less than $1\times10^{-13}$ for the whole domain.

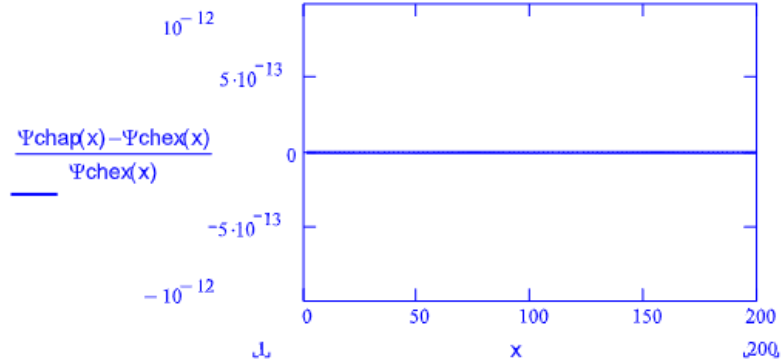

**Fig. 26** Plot of the relative error of our approximant to the Chebychev function, w.r.t. the exact function, in the domain *[0,200]* showing the high accuracy, order of $10^{-7}$, of our approximant.

In Fig. A-2 of the Appendix 2 we again compare both functions, Eqs (44) and (45), but for a domain of higher values of the independent variable, namely *x∈[5000, 5050]*, those plots showing the same order of accuracy of our approximant in Eq (45).

**7. Riemann Hypothesis and the Chesbyshev Function of Second Kind**

The main open question about prime numbers is their relation to the Riemann Hypothesis. A well-known theorem that connects Chebyshev's second class function to the Riemann Hypothesis was presented by L. Shoenfield [7] in 1976. He then proved that:

"The Riemann hypothesis implies the inequality:

$$"|\Psi_{Che}(x) - x| < \frac{1}{8\pi}\sqrt{x}\, log^2(x)", \quad \text{if } x \geq 73.2\," \qquad (46)$$

This is in fact an important sharp inequality that tells us: Should the Riemann Hypothesis holds then we could prove the inequality estimate.  In Fig. 27 we have plotted this inequality using both, our approximant to the Chebyshev function (on the left), which is based on our prime numbers generating function, and for comparison using the actual exact Chebyshev function (on the right), both plots in the real domain (0, 400).

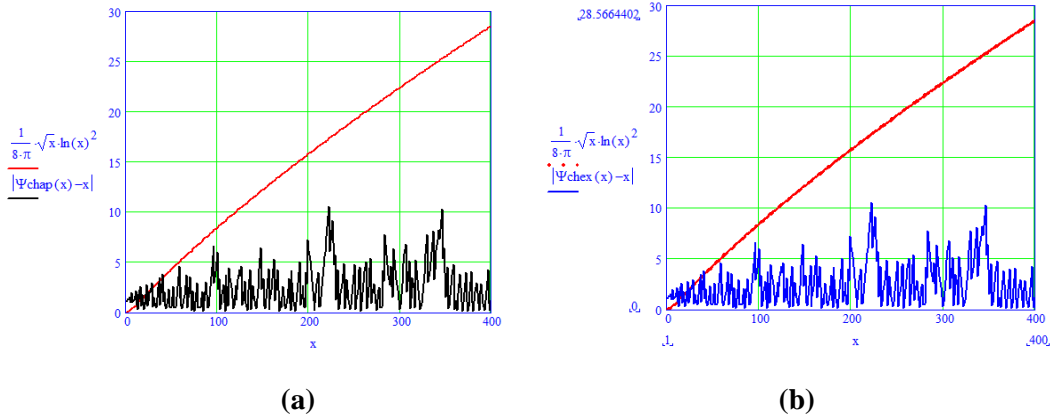

**(a)** **(b)**

**Fig. 27** Plots *vs. x*∈*R* of Shoenfield`s sharp inequality in Eq. (46). **(a)** Using our approximant to Chebyshev 2[nd] class function; and **(b)** Using the exact expression for that function. The two plots are indistinguishable. The inclined line plotted in both figures represents the function $\frac{1}{8\pi}\sqrt{x}\, ln^2(x)$ that belongs on the right side of the theorem inequality.

It may be seen that both figures are coincident, and that the inequality holds in both cases: note that for x< 73 the oscillating curve intersects the line above it, while for x>73 the oscillating function is always below the line that represents the function $\frac{1}{8\pi}\sqrt{x}\, ln^2(x)$, verifying that Riemann Hypothesis must hold.

## 8. Results and Discussion

In this work we derived a new Prime numbers Generating Function that gives the exact number of primes in any given integer domain. This generating function of primes is defined using a new prime number discriminating function Λ, the latter obtained in terms of elementary functions, presented and duly derived in Section 2. Not dependent upon prime numbers tables, or algorithms, these two functions are easier-to-use tools for further prime numbers arithmetical function and relations, as well as for applications in basic sciences e.g. in physics, as in the applications presented in Section 5. We presented solid evidence (see Table 1, Appendix 1) of the exactitude of our prime number generating function by generating 64 random integers lying in the arbitrarily chosen integer domain *[13, 2591491]*, and then successfully applying our prime discriminating function *Λ* to determine whether these random integers belong or not to the primes distribution in that interval (Section 2, Table 1 in Appendix 1). Also introduced was a new Prime Numbers Counting Function *C(2,u)*, giving the exact number of primes from 2 and up the integer *u*, that can be easily implemented using a personal computer. This counting function happens to be exact, generating the well-known staircase-like plot of primes when plotted for any integer domain. Both, our new prime numbers counting and generating functions have been used here to derive finite approximants to Euler's Product formula and to the Chebyschev function of second kind [Eq. (51)].These approximants are shown to be of good accuracy, achieving relative errors of order $10^{-7}$, even when calculated with relatively small upper bound values (*e.g. H=100* in Eq.(*25*)) of the finite product or sum in their definitions, and using modest personal computers. Much better results should be obtained using mainframe computers nowadays available at large, to run the computer programs implementing our three Prime Numbers approximants with a high upper bound, say *H=$10^6$*, or even higher. We also proceeded to generate an accurate, finite-sum approximant to Riemann Zeta function $\zeta(\sigma,\tau)$, using our prime numbers discriminating function, that allowed us to find the location of non-trivial zeroes of the Zeta function along the critical line *($\sigma_0$=0.5, $\tau$ )* with great accuracy, for both the location of the non-trivial roots −with accuracies of the order of 0.02 % − and for the expected zero values of the modulus of the Zeta function at those non-trivial roots, with accuracies in the range of $10^{-10}$, or even better. In Section 5 we applied the Principle of Least Action to obtain the coordinates *($\sigma,\tau$)* of non-trivial roots of function Zeta on the complex plane, and to the effect we considered a simple harmonic oscillator system, and a variational test function based on the product of the *squared modulus o*f our approximated Zeta function times *sin (t)*, just to fulfil border conditions of the required interval of definition. With these we obtained the required Lagrangian (Eq (34)) to then apply the Principle of Least Action, which rendered us an exact analytical and simple Action function, with which we then found the sought coordinates of non-trivial roots of the Zeta function with good accuracy. The minimum values of our obtained Action at these roots are close to zero within $10^{-5}$ or even less, as we expected. We then extended our application of the Least Action principle by defining, and using, a variational test function written again in terms of the modulus of the Riemann zeta function, yet this time not squared, but to the power $1/\sigma$. The new results showed that the single compatible value of the coordinate σ within our formalism is *$\sigma_0$=0.5*, which happens to be the particular value precisely required by the Riemann Hypothesis to hold. We also derived an important symmetrical Functional Equation that is, in good sense, an analogous equation to the well-known Riemann Functional Equation. This new function allowed us to obtain in a simple way non-trivial roots of our approximant to Riemann zeta function.

Finally, we defined and used a finite self-consistent approximant to the Chebyschev function of 2$^{nd}$. kind, that does not require prime number tables or numerical algorithms, but that instead is written in terms of our prime number counting function *C(2,x)*, to successfully replace the the $\pi(x)$ counting and the Chebyschev function, in the well-known Schoenfield's and Trudgian's theorems relations [ 5, 6, 7] that prove sharp inequalities involving those two functions, and that once satisfied imply that Riemann Hypothesis holds.

## 9. Acknowledgements

We acknowledge Lcdo. Frederick A. Luy Alluera of the Telematics Division of Universidad Simón Bolívar (USB), for his assistance to improve our connectivity to Internet. We also thank our Physics Majors students and our Physics Bachelors, the latter now doing graduate work abroad, who organized in *cofAlumnae* of Universidad Simón Bolívar (USB), have provided support to improve our connectivity, particularly to Miss. Karleyda Sandoval (B. Sc. Physics, USB) for managing that support.

**Appendix 1**

The randomly-generated Table I of integers shows the exactitude of our Prime Number Generating function $\psi$ introduced in Section 2. Four sets of 16 *random integers* $u_n$, with $n \in [1,64]$ are listed in the four 2$^{nd}$ columns of the sets: 16 integers per column. These 64 integers $u_n$ were generated using a computer-programmed discrete differential random formalism (described in the following paragraph). The program uses as input a pair of integers parameter $(K_i, u_i)$ to generate each set of 16 random integers). The integers $\psi(u_n)$ returned by our prime generation function $\psi$, for each random integer $u_n$ in the table, appear in the third column of each set: if $u_n$ is prime our function simply returns that $u_n$, otherwise it returns a zero, e.g. for the $n^{th}=16^{th}$ run of the program the random formalism gave us the integer $u_{n=}1270197$ which is not prime (since it factorizes as $1319 \times 107 \times 9$) and you find the expected $\psi(1270197)=0$ in Table 1.

## TABLE 1

| n = | $u_n$ = | $\Psi(u_n)$ = | n = | $u_n$ = | $\Psi(u_n)$ = |
|---|---|---|---|---|---|
| 1 | 35 | 0 | 17 | 13 | 13 |
| 2 | 71 | 71 | 18 | 27 | 0 |
| 3 | 143 | 0 | 19 | 57 | 0 |
| 4 | 289 | 0 | 20 | 115 | 0 |
| 5 | 583 | 0 | 21 | 241 | 241 |
| 6 | 1169 | 0 | 22 | 491 | 491 |
| 7 | 2351 | 2351 | 23 | 1023 | 0 |
| 8 | 4735 | 0 | 24 | 2139 | 0 |
| 9 | 9565 | 0 | 25 | 4311 | 0 |
| 10 | 19173 | 0 | 26 | 9071 | 0 |
| 11 | 38597 | 0 | 27 | 18729 | 0 |
| 12 | 77915 | 0 | 28 | 38455 | 0 |
| 13 | 156137 | 0 | 29 | 81473 | 0 |
| 14 | 312863 | 312863 | 30 | 168543 | 0 |
| 15 | 631507 | 631507 | 31 | 347025 | 0 |
| 16 | 1270197 | 0 | 32 | 730789 | 730789 |

| n = | $u_n$ = | $\Psi(u_n)$ = | n = | $u_n$ = | $\Psi(u_n)$ = |
|---|---|---|---|---|---|
| 33 | 73 | 73 | 49 | 9 | 0 |
| 34 | 147 | 0 | 50 | 19 | 19 |
| 35 | 297 | 0 | 51 | 39 | 0 |
| 36 | 597 | 0 | 52 | 81 | 0 |
| 37 | 1195 | 0 | 53 | 173 | 173 |
| 38 | 2403 | 0 | 54 | 361 | 0 |
| 39 | 4831 | 4831 | 55 | 765 | 0 |
| 40 | 9745 | 0 | 56 | 1661 | 0 |
| 41 | 19619 | 0 | 57 | 3529 | 3529 |
| 42 | 39571 | 0 | 58 | 7457 | 7457 |
| 43 | 79229 | 79229 | 59 | 15189 | 0 |
| 44 | 158961 | 0 | 60 | 32065 | 0 |
| 45 | 318841 | 318841 | 61 | 65687 | 65687 |
| 46 | 638585 | 0 | 62 | 139319 | 0 |
| 47 | 1287937 | 0 | 63 | 294927 | 0 |
| 48 | 2591491 | 0 | 64 | 619019 | 619019 |

**Integer number generator using a random differential discrete formalism**

To generate Table 1 the program listed below was applied in the integer domains *[1, 16 ], [17,32 ], [33,48] and [49,64 ]*. For each domain the program first generates a pair of random integers $K$ and $u$ (shown in Table 1 as $u_1$, $u_{17}$, $u_{33}$ and $u_{49}$)

$$n = 1, 2, 3, \ldots 16$$

$$K = 1 + trunc[rnd(100)]$$

$$u_1 = 2\, trunc\left[rnd\left(\frac{K}{2}\right)\right] + 1$$

$$u_{n+1} = \left\{2\left[u_n + trunc\left(rdn\left(\frac{u_n}{K}\right)\right)\right] + 1\right\}$$

**Appendix 2**

Below, and for comparison, we plot the Chesbyshev exact function, Eq (44) and our approximant Eq (45) to that function in the real domain [5000, 5050]

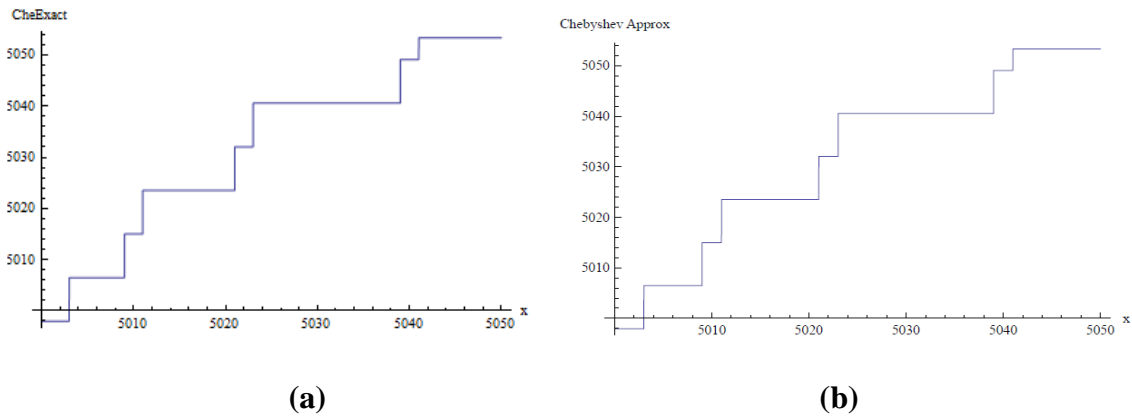

(a)  (b)

**Fig. A2** Comparison plots of the exact Chebyschev function **(a)** in Eq. (45), and of our approximant **(b),** in Eq. (46), in the domain *[5000, 5050]* with relative accuracy order of $10^{-7}$.

**Appendix 3**

We may use our procedure, based on the Hamilton Principle, and presented in Section 5, to find non-trivial roots of the Zeta function in a short domain in the neighbourhood of the higher ordinate $\tau = 999$, domain where we know. from accepted tables of non-trivial roots of the Riemann zeta function, there exists one of such roots. We need to plot the action $A(\sigma, \tau)$ versus the pair of variables: $\omega = \tau - \sigma$ and $\eta = \tau + \sigma$ defined in the main text (Sub-section 5.1). To improve the accuracy of our calculation we now use the upper bound $n_{max} = 1000$ in the Lagrangian of Eq. (33). In fact we zoom in a short domain about the ordinate $\tau = 1000$. The obtained plots of the action are presented in Fig. A3(a) and A3(b), respectively.

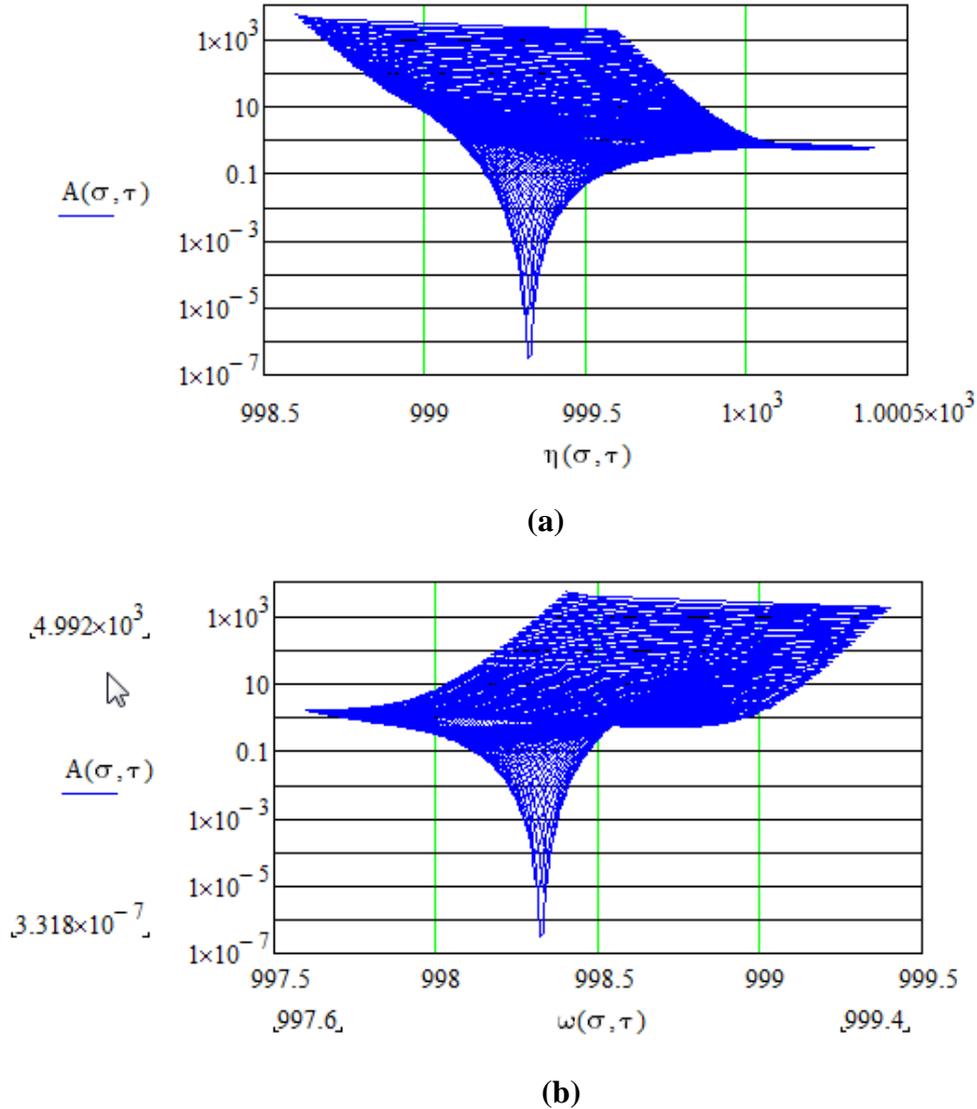

**Fig. A3** Plots of the action $A(\sigma, \tau)$ versus the variables $\eta = \tau + \sigma$ and $\omega = \tau - \sigma$, sum and difference of the real and imaginary parts of the complex $s=(\sigma, \tau)$ in figures **(a)** and **(b)**, respectively.

In Fig. **A3(a)** the minimum of the Action occurs at $\eta = 998.32$, while in Fig. **A3 (b)** the minimum occurs at $\omega = 999.32$. Thus we may write:

$$\tau = \tau - \sigma = 998.32 \quad \text{and} \quad \omega = \tau + \sigma = 999.32.$$

These equations give us the two coordinate's values: $\sigma = 0.495$, $\tau = 998.825$. According to the Riemann Hypothesis the expected correct values are: $\sigma = 0.5$ and $\tau = 998.827547137$. This shows the good accuracy of our method based on analytical mechanics, in fact the calculated relative errors ($\varepsilon$) of our results are rather small:

$$\varepsilon_\sigma = 0.010000 \quad \text{and} \quad \varepsilon_\tau = 0.0000026$$